\documentclass{article}
\usepackage{amsfonts,amssymb,amsmath} 
\newtheorem{Lemma}{Lemma}[section]\newcommand{\bel}{\begin{Lemma}}\newcommand{\eel}{\end{Lemma}}
\newtheorem{Proposition}[Lemma]{Proposition}\newcommand{\bprop}{\begin{Proposition}}\newcommand{\eprop}{\end{Proposition}}
\newtheorem{Theorem}[Lemma]{Theorem}\newcommand{\bthe}{\begin{Theorem}}\newcommand{\ethe}{\end{Theorem}}
\newtheorem{Remark}[Lemma]{Remark}\newcommand{\beR}{\begin{Remark}\rm}\newcommand{\eeR}{\end{Remark}}
\newtheorem{Definition}[Lemma]{Definition}\newcommand{\bed}{\begin{Definition}}\newcommand{\eed}{\end{Definition}}
\newtheorem{Example}[Lemma]{Example}\newcommand{\bex}{\begin{Example}\rm}\newcommand{\eex}{\end{Example}}
\newtheorem{Corollary}[Lemma]{Corollary}\newcommand{\bcor}{\begin{Corollary}\rm}\newcommand{\ecor}{\end{Corollary}}

\newcommand{\bpr}{{\bf proof:~}}\newcommand{\epr}{\hfill$\blacksquare$\\}  

\newcommand{\beq}{\begin{equation}}\newcommand{\eeq}{\end{equation}}
\newcommand{\bem}{\begin{displaymath}}\newcommand{\eem}{\end{displaymath}}
\newcommand{\beqa}{\begin{eqnarray}}\newcommand{\eeqa}{\end{eqnarray}}
\newcommand{\bee}{\begin{enumerate}}\newcommand{\eee}{\end{enumerate}}
\newcommand{\bei}{\begin{itemize}}\newcommand{\eei}{\end{itemize}}
\newcommand{\ber}{\begin{array}{l}}\newcommand{\eer}{\end{array}}
\newcommand{\bet}{\begin{tabular}{l}}\newcommand{\eet}{\end{tabular}}
\newcommand{\bpm}{\begin{pmatrix}}\newcommand{\epm}{\end{pmatrix}}

\newcommand{\di}{\partial}
\newcommand{\li}{~\\ $\bullet$ }

\newcommand{\tinyM}{\scriptstyle}
\newcommand{\tinyT}{\scriptsize}

\newcommand{\ra}{\rightarrow}

\newcommand{\cU}{{\cal{U}}}
\newcommand{\cV}{{\cal{V}}}

\newcommand{\he}{\hat{e}}

\newcommand{\mC}{\mathbb{C}}\newcommand{\mD}{\mathbb{D}}\newcommand{\mL}{\mathbb{L}}\newcommand{\mN}{\mathbb{N}}
\newcommand{\mP}{\mathbb{P}}\newcommand{\mPN}{\mP_f^{N_d}}\newcommand{\mR}{\mathbb{R}}\newcommand{\mS}{\mathbb{S}}
\newcommand{\mV}{\mathbb{V}}\newcommand{\mZ}{\mathbb{Z}}

\newcommand{\tSi}{\widetilde{\Sigma}} \newcommand{\tV}{\widetilde{V}}\newcommand{\tcV}{\widetilde{\cV}}

\newcommand{\al}{\alpha}
\newcommand{\si}{\sigma}
\newcommand{\Si}{\Sigma}

\newcommand{\gNnd}{generalized Newton-non-degenerate } \newcommand{\Nnd}{Newton-non-degenerate }

\newcommand{\ellipse}[8]
{\put(0,0){\qbezier(#1,#2)(#1,#4)(#3,#4)}\put(0,0){\qbezier(#3,#4)(#5,#4)(#5,#6)}\put(0,0){\qbezier(#5,#6)(#5,#8)(#7,#8)}
\put(0,0){\qbezier(#7,#8)(#1,#8)(#1,#2)}
}

\title{Enumeration of uni-singular algebraic hypersurfaces
}
\author{D. Kerner
}

\begin{document}
\setcounter{secnumdepth}{6} \setcounter{tocdepth}{2}
\newcounter{temp}\newcounter{kk}\newcounter{iter}\newcounter{Vh}\newcounter{Ch}
\newcounter{tempx}\newcounter{tempy}\newcounter{Lx}\newcounter{Ly}\newcounter{Lslope}
\newcounter{tx1}\newcounter{tx2}\newcounter{tx3}\newcounter{ty1}\newcounter{ty2}\newcounter{ty3}
\newcounter{x1}\newcounter{x2}\newcounter{x3}\newcounter{y1}\newcounter{y2}\newcounter{y3}

\maketitle
\begin{abstract}
We enumerate complex algebraic hypersurfaces in $\mP^n$, of a given (high) degree with one singular point of
a given singularity type. Our approach is to compute the (co)homology classes of the corresponding equi-singular
strata in the parameter space of hypersurfaces. We suggest an inductive procedure, based on intersection theory
combined with liftings and degenerations. The procedure computes the (co)homology class in question, whenever a
given singularity type is properly defined and the stratum possesses good geometric properties. We consider in
details the generalized Newton-non-degenerate singularities. We give also examples of enumeration in some other cases.
\end{abstract}
\tableofcontents
\section{Introduction}
\subsubsection{Preface}
This paper is a generalization of the previous one \cite{Ker}, where enumeration of plane singular curves (with
one singular point) was done.
The results for curves were as follows:
\li for a large class of singularity types (the so-called linear singularities) we gave a method to write immediately
 explicit formulae solving enumeration problem.
\li for all other singularity types we gave an algorithm (which
is quite efficient and for every particular singularity type gives the final answer in a bounded number of steps).

The goal of this work is to generalize the method to the case of (uni-singular, complex, algebraic) hypersurfaces
in $\mP^n$.

The theory of singular hypersurfaces is much richer and complicated than that of curves. Correspondingly the
enumeration is much more difficult both technically and conceptually.
It seems that there does not exist one (relatively) easy method applicable to all types of singularities.
We propose a method of calculation applicable to a class of the {\it generalized Newton-non-degenerate}
singularities (this includes in particular the $A,D,E$ types and all the singularities with $\mu\le14$).
The method can also be applied to some other singularity types, we consider examples in Appendix A.
\subsubsection{General settings}
We work with (complex) algebraic hypersurfaces in the ambient space $\mP^n$. A hypersurface
is defined by a polynomial equation $f(x)=0$ of degree $d$ in the homogeneous coordinates of $\mP^n$.
The parameter space of such hypersurfaces (the space of homogeneous polynomials of total degree $d$
in $(n+1)$ variables) is a projective space.
We denote it by $\mPN$ (here \mbox{$N_d={d+n\choose{n}}-1$} and $d$ is assumed to be sufficiently high).

The {\it discriminant}, $\Sigma\subset\mPN$, is the (projective) subvariety of the parameter space,
whose points correspond to the singular hypersurfaces (generic points of the discriminant
correspond to hypersurfaces with one node). Everywhere in this paper (except for Appendix A)
we restrict consideration {\it to isolated} singularities. Even more,
we consider only hypersurfaces with just {\it one} singular point.

When working with a singular point (specifying its type, parameters etc.) we usually pass from the
category of projective hypersurfaces to that of hypersurface germs and consider everything in a
small neighborhood in classical topology.

Consider the classification by (local embedded) topological type: two hypersurface germs $(V_i,0)\subset(\mC^n,0)$
are of the same type if there
exists a homeomorphism $\mC^n\stackrel{\phi}{\ra}\mC^n$ such that $\phi(V_1)=V_2$.
Two singular projective hypersurfaces are said to be of the same (local embedded topological) singularity type $\mS$,
if the corresponding germs are.
For a given topological type $\mS$, consider the stratum of (projective) hypersurfaces $\Si_\mS$ with
a singular point of this type.

Unlike the case of curves the notions of the topological equivalence and the corresponding equisingular
stratification are quite complicated for hypersurfaces (as is discussed shortly in $\S$ \ref{SecOnEquisingularity}).
For the purposes of enumeration we use a more restricted equivalence: by the Newton diagram.
We work always with {\it commode} (convenient) diagrams.

A singular hypersurface germ is called {\it generalized Newton-non-degenerate} if it can be brought
to a Newton-non-degenerate form by locally analytic transformations.
A hypersurface with one singular point is called \gNnd if the corresponding germ is \gNnd.
 (For the precise definitions and discussion of the
relevant notions from singularities cf. $\S$ \ref{SecTypesSingul}).
Everywhere in this paper (except for Appendix A) we consider \gNnd hypersurfaces.

\bed\label{DefNDEquisingularStratum}
Two (generalized Newton-non-degenerate) hypersurface-germs are called ND-equivalent if they can be brought by locally
analytic transformations to \Nnd forms with the same Newton diagram.

For a given Newton diagram $\mD$, the equisingular family $\Sigma_{\mD}^{d,n}$ is defined as the
set of all the points in the parameter space $\mPN$, corresponding to
\gNnd hypersurfaces of degree $d$ that can be brought by locally analytic transformations to $\mD$.
\eed
Note that for \gNnd hypersurfaces this equivalence implies equivalence by the embedded
topological type. In general this equivalence is weaker than the (contact) analytic equivalence.
Therefore we call this equivalence (and the corresponding families) {\it ND-topological}.
From now on (except for $\S$ \ref{SecOnEquisingularity}), by singularity type we mean the
ND-topological type $\mD$.

As the degree of the hypersurfaces $d$ and the dimension $n$ are always fixed we omit them.
 For the enumeration purposes we always work with the {\it topological closure} of the strata
 ($\bar\Si_\mD\subset\mPN$), to simplify the formulae we usually omit the closure sign
 (e.g. $\Sigma_{A_k}$, $\Sigma_{D_k}$, $\Sigma_{E_k}$ etc.).

The so defined closures sometimes coincide with the closures of the topological equisingular
strata (section \ref{SecOnEquisingularity}).
For example this is the case for $A,D,E$ singularities \cite{AVGL,Var}.

To specify the ND-topological singularity type (i.e. to construct the corresponding diagram) we usually give
 a representative. In the simplest cases this representative (the normal form) is classically fixed
 (cf. tables in \cite[chapter 1]{AVGL}).
\bex
The normal forms of some simplest singularities (since we do not consider analytical equivalence
the moduli are omitted):
\beq\ber\label{NormalForms}
A_k:~z_1^{k+1},~~D_k:~z_1^2z_2+z_2^{k-1},~~E_{6k}:~z_1^3+z_2^{3k+1},~~E_{6k+1}:~z_1^3+z_1z_2^{2k+1},~~E_{6k+2}:~z_1^3+z_2^{3k+2}
\\
P_8:~z_1^3+z_2^3+z_3^3,~~X_9:~z_1^4+z_2^4,~~J_{10}:~z_1^3+z_2^6,~~T_{p,q,r}:~z_1^p+z_2^q+z_3^r+z_1z_2z_3,~~
\frac{1}{p}+\frac{1}{q}+\frac{1}{r}<1
\\
Q_{10}:~z^3_1+z^4_2+z_2z^2_3,~~S_{11}:~z^4_1+z^2_2z_3+z_1z^2_3,~~U_{12}:~z^3_1+z^3_2+z^4_3
\eer\eeq
Here we consider the normal form up to stable equivalence, i.e. up to the (non~degenerate) quadratic forms:
$f(z_1,\dots,z_k)+\sum_{i=k+1}^nz_i^2\sim f(z_1,\dots,z_k)$.
\eex

Everywhere in the paper we assume that the degree of the hypersurfaces, $d$, is high. A
sufficient condition is: $d$ is bigger than the degree of determinacy for a
given singularity type (the later is often the maximal coordinate of the Newton diagram).
This condition is not necessary, e.g. in the case of curves the method works also in the
irregular region (of small $d$), the algorithm should be slightly modified \cite[section 5]{Ker}.

The so defined equsingular strata are (almost by construction): non-empty, algebraic, pure dimensional
 and irreducible (cf. proposition \ref{ClaimNDtopStratumIsAlgVariety}).
 Therefore the enumerative problem is well defined.

Every equisingular stratum is (by construction)
embedded into $\mPN$, we identify the stratum with its embedding. Its closure has the homology class in
 the corresponding integer homology group:
\beq
[\overline{\Sigma}_\mD]\in H_{2dim(\Si_\mD)}(\mPN,\mZ)\approx\mZ
\eeq
The degree of this class is the degree of the stratum.
The {\it cohomology class of the stratum}, $[\overline{\Sigma}_\mD]\in H^{2N_d-2dim(\Si_\mD)}(\mPN,\mZ)$, is the class dual (by Poincare
duality) to the above homology class.

We denote the homology and cohomology class by the same letter, no confusion should arise.
\subsubsection{The goal of the paper, motivation and main results}
The goal of this paper is to provide a method to enumerate the hypersurfaces, i.e. to calculate the cohomology
classes, for the strata corresponding
to hypersurfaces with one singular point of a given ND-topological type (for generalized
Newton-non-degenerate hypersurface-germs).
\\
\\
\\
The discriminant and more generally, varieties of equisingular hypersurfaces (Severi-type varieties),
have been a subject of study for a long time.
Already in 19'th century it was known, that the (closure of the) variety of nodal hypersurfaces of degree d in $\mP^n$
(i.e. the discriminant)
is an irreducible subvariety of $\mPN$ of codimension 1 and degree:
\beq
(n+1)(d-1)^n
\eeq
Any further progress happens to be difficult. The work was mainly concentrated on the enumeration of curves on surfaces
\cite{KleiPien1} with many simple singularities.

The present situation in enumeration of singular hypersurfaces seems to be as follows.
(This is not a complete/historical review, for a much better description cf: \cite{Klei1,Klei2,Kaz4}.)
\li{In 1998 P.Aluffi\cite{Alufi}} has calculated the degrees of the strata of cuspidal and bi-nodal hypersurfaces
 ($\Si_{A_2}$, $\Si_{2A_1}$)
\li{In 2001 R.Hern\'aández and M.J.V\'aázquez-Gallo\cite{HV}} enumerated most of the singularities
 of cubic surfaces in $\mP^3$.
\li{In 2003 I.Vainsencher\cite{Vain2}} has calculated the degrees of some strata of multi-nodal hypersurfaces
($\Si_{rA_1}$ for $r\leq6$).
\li{In 2000-2003 M.Kazarian in the series of papers \cite{Kaz1,Kaz2,Kaz3,Kaz4}} used topological approach to
prove that there
exist {\it universal} formulae for the degrees of equisingular strata. In the spirit of Thom \cite{Thom54},
they are (unknown) polynomials in some combinations of the Chern classes of the ambient space
and the linear family. (In our case these are the Chern classes of $\mP^n$, and $\mPN$.)
The coefficients of those polynomials depend on the singularity type only (and not on the degree of
hypersurfaces or their dimension). The enumerative answer for a particular question is obtained just by
substitution of the corresponding Chern classes into the universal polynomial.

Kazarian has developed a method for calculation of these Thom polynomials. From them one gets the
degrees for strata of singular hypersurfaces. In particular, he gives the answers for all the possible
combinations of types up to codimension 7.

His method enumerates all the singularity types of the given codimension simultaneously. Therefore
it needs a preliminary classification of the singularities of a given codimension. Even if one does this,
the computations are non-effective when one needs the answer for just one type (e.g. $A_k$)

Our motivation was to show that the approach suggested in \cite{Ker} (and used there to completely solve the problem
for plane uni-singular curves) generalizes to the case of hypersurfaces. In particular, from the theorems
\ref{TheoremEnumLinearSings},\ref{TheoremEnumNonLinearSings} it follows
\bcor\label{ClaimMainOfTheMethod}
The proposed method of degenerations (the algorithm) allows
enumeration of any (generalized Newton non-degenerate) singularity (in a bounded number of steps).
\ecor

The result of the enumeration procedure is the cohomology class of a stratum $Si_\mD$: a polynomial of degree $n$
(the dimension of the ambient space, $\mP^n$) in $d$ (the degree of hypersurfaces).
The coefficients are functions of $n$ and of the singularity type.
As an example of calculations we have
\bprop
The cohomology classes of the lifted strata in the following cases are given in Appendix \ref{SecApendCohomologClasses}:
 ordinary multiple points, reducible multiple points (see the definition
in Appendix \ref{SecReducibleForms}), $A_{k\leq4}$, $D_{k\leq6}$, $E_{6}$, $P_8, X_9$, $Q_{10}$, $S_{11}$,
$U_{12}$
\eprop
From our results one can obtain some restrictions on the universal Thom polynomials. We should note, however,
that from our results it is impossible to recover Thom polynomials completely when the $\mu\ge7$ (cf. Appendix C).
\subsubsection{Description of the method}
Here we briefly describe the method. It is considered in more details in $\S$ \ref{SecMainIssues}.
We start from the ingredients and then formulate the enumeration theorems.

We go in a naive way, trying to work with (locally) complete intersections of  hypersurfaces
defined by explicit equations. The resulting cohomology class is obtained as the product of the classes of
hypersurfaces, with various corrections subtracted. Repeat that we always work with the closed strata.
\paragraph{\hspace{-0.3cm}Liftings.\hspace{-0.3cm}}\label{SecMinimalLifting}
To write explicit equations,
 lift a given equisingular stratum (which initially lies in $\mPN$) to a bigger space ($Aux\times\mPN$). Here $Aux$ is an
auxilliary space that
traces the parameters of the singularity (singular point, tangent cone etc.).
\bex The {\it minimal lifting (partial desingularization)} is just the universal hypersurface
\beq\label{ExampleLiftings}
\tSi_\mD(x):=\overline{\left\{(x,f)\Big|\ber \mbox{The~hypersurface~defined~by~}f(x)=0~\mbox{has}\\
\mbox{singularity type}~\mD~\mbox{at~the point~}x\eer\right\}}\subset\mP^n_x\times\mPN
\eeq
Here $\mP^n_x$ is the ambient space of singular hypersurfaces (the subscript $x$ emphasizes that the point of the space
 is denoted by $x$).
\eex
The cohomology class of the lifted version ($\tSi_\mD$) is easier to calculate (e.g. for ordinary
multiple point the lifted stratum is just a complete intersection, cf. $\S$ \ref{SecOrdinaryMultiplePoint}).
The (co)homology class is now not just a number, but a polynomial (in the generators of the cohomology ring
of the bigger ambient space). Hence we have the multidegree of $\tSi$.
This provides, of course, much more information about a particular stratum.
\\\\
Once the class $[\tSi_\mD]$ has been calculated, the cohomology class of the
 original stratum ($[\Si_\mD]$) is obtained using the Gysin homomorphism. Namely, the projection $Aux\times\mPN\stackrel{\pi}{\rightarrow}\mPN$
induces the projection on homology: $H_{i}(Aux\times\mPN)\stackrel{\pi_*}{\rightarrow}H_i(\mPN)$. By Poincare duality this gives a projection
in cohomology: $H^{i+2dim(Aux)}(Aux\times\mPN)\stackrel{\pi_*}{\rightarrow}H^i(\mPN)$.
It sends the component $H^k(\mPN)\otimes H^{2dim(Aux)}(Aux)$ isomorphically to $H^k(\mPN)$ and sends all other cohomology classes to zero.

From the calculational point of view, we should just extract from the cohomology class $[\tSi_\mD]$
(which is a polynomial in the cohomology ring of
$Aux\times\mPN$) the coefficient of the maximal powers of the generators of $H^*(Aux)$.

In the above example of minimal lifting this homomorphism is $H^{k+2n}(\mP^n_x\times\mPN)\mapsto H^{k}(\mPN)$,
and from $[\tSi_\mD]$ one should extract
the coefficient of the $n$'th power of the generator of $H^*(\mP^n_x)$.

Summarizing: the cohomology class of a stratum $\Sigma_\mD$ is completely fixed by that of its lifted version $\tSi_\mD$.
\paragraph{\hspace{-0.3cm}Degenerations from simple to complicated.\hspace{-0.3cm}}The lifted stratum can often be globally defined
by some explicit equations (the case of {\it linear singularity}, cf. definition \ref{DefLinearSingularities}
 in $\S$ \ref{SecLinearSingularDefinitions}).
 Unfortunately it is usually only a locally complete intersection (but not globally).
So, if one chooses a locally defining set of hypersurfaces, their intersection contains some residual (unnecessary)
pieces, whose contribution to the cohomology class should be subtracted. The serious difficulty is that these
residual pieces can be of dimension {\it bigger} than the true stratum.

In this case we proceed as follows.
Let $\mD$ be the singularity under consideration, let $\mD_0$ be some singularity type to which $\mD$ is adjacent
and for which the enumeration is already done (the trivial choice for $\mD_0$ is just an ordinary multiple point of
the same multiplicity as $\mD$). Represent $\mD$ as a chain of successive degenerations (i.e. adjacencies),
starting from $\mD_0$.
At each step the codimension of the variety grows by 1, the stratum being intersected by a hypersurface.
Each intersection can be non-transversal somewhere, the resulting variety of the intersection is usually reducible.
In addition to the needed (true) variety it contains some residual varieties. We emphasize, that at
each step the intersection is with a hypersurface, therefore the dimensions of residual pieces are
not bigger than that of the true variety.
Thus, at each step the contribution of residual pieces can be removed from the cohomology class of the intersection.
In this process one should check:
\li  where the non-transversality happens? This question is considered in $\S$ \ref{SecPossibleCyclesOfJump}
using proposition \ref{ClaimTransversalityUnderNewtonDiagram}.
\li what are the residual pieces produced in the intersection? This question is considered in sections
\ref{SecLiftings} and \ref{SecIdeologyOfDegenerations}.
\li how "to remove" their contributions from the answer? This amounts to calculation of the cohomology classes of
residual varieties and is considered in sections \ref{SecCohomologyClassesOfCyclesofJump}
and \ref{SecCohomologyClassOfRestrictionFibration}.

All above can be formulated as a proposition (proved in \ref{SecIdeologyOfDegenerations}):
\bthe\label{TheoremEnumLinearSings}
\li
For a given linear type $\mD$ and an auxiliary type $\mD_0$ the algorithm forms the chain of degenerations
$\mD_0\ra..\ra\mD$.
All the vertices correspond to linear types (fixed by the choice of $\mD,\mD_0$).
The number of vertices equals $codim(\Si_{\mD})-codim(\Si_{\mD_0})$.
\li Each edge $\mD_i\ra\mD_{i+1}$ (a codimension 1 degeneration) provides a linear expression
 for the class $[\mD_{i+1}]$ in terms of $[\mD_i]$ and the classes of residual cycles. The residual cycles are
 fixed by the geometry of the degeneration.
\li The number of steps, needed to achieve the result, is not bigger than the number of points under the Newton diagram.
\ethe
\paragraph{\hspace{-0.3cm}Degenerations from complicated to simple (simplifying degenerations).\hspace{-0.3cm}}\label{SecSimplifyingDegenerations}
 In most cases, even the lifted stratum is difficult to define explicitly. This is the case of
{\it non-linear} singularities (defined in $\S$ \ref{SecLinearSingularDefinitions}). Then, instead of trying to arrive at the needed
singularity $\mD$ by degenerations of some
 simpler singularity, we degenerate the $\mD$ itself. The goal is to arrive at some singularity ($\mD'$) of higher
 codimension (or higher Milnor number),
 which is however simple to work with (e.g. linear singularity).

Or, geometrically, we intersect the lifted stratum ($\tSi_{\mD}$) with a cycle
in the ambient space so that the cohomology
class of the resulting stratum ($\tSi_{\mD'}$) is easier to calculate. Then (if the intersection is transversal)
we have the equation for the cohomology classes:
\beq
[\tSi_{\mD}]\times[\mbox{degenerating cycle}]=[\tSi_{\mD'}]\in H^*(Aux\times\mPN)
\eeq
We choose the degenerating cycle in the appropriate manner, so that the above equation of cohomology classes
fixes the class of $\tSi_{\mD}$ uniquely (cf. $\S$ \ref{SecInvertibilityofDegeneration}).

In fact the situation is more complicated:
\li The resulting stratum in general is reducible and non reduced. Its reduced components enter with different
multiplicities (since the initial stratum $\tSi_{\mD}$ is singular at these loci). The resulting strata and
their multiplicities are obtained from the defining ideal, by explicit check of the equations.
\li The needed stratum $\tSi_{\mD'}$ usually is not a complete intersection.
Thus, on the right hand side in the above equation there can appear some residual pieces. In this case one should also
remove their contributions.

The result for non-linear (\gNnd) singularities is
(proved in $\S$ \ref{SecIdeologyOfDegenerations}):
\bthe\label{TheoremEnumNonLinearSings}
\li For each (non-linear) singularity type $\mD$ a tree of degenerations is constructed.
The root of the tree is the original type $\mD$, the
leaves are some linear singularities, adjacent to the original stratum. This tree is constructed from the
Newton diagram of the
given singularity type, without any preliminary classification or preliminary knowledge of adjacent strata.
\li Every edge of the tree corresponds to a degeneration, resulting in a pure dimensional variety.
The corresponding equation for cohomology
classes is of the form  $[\tSi_{\mD_i}]\times[\bet degenerating\\divisor\eet]$=$\sum a_j[\tSi_{\mD_{i+1,j}}]+[\bet residual\\piece\eet]$.
The cohomology class of the residual piece is calculated by a standard procedure. The cohomology class of the stratum
$[\tSi_{\mD_i}]$ is
restored uniquely from this equation.
\li
If the initial non-linear singularity has order of determinacy
$k$ and multiplicity $p$, then the number of vertices in this tree is not bigger than ${k+n\choose{n}}-{p-1+n\choose{n}}$.
\ethe
\paragraph{\hspace{-0.3cm}Some special simple cases.\hspace{-0.3cm}}
In some (very special) cases the enumeration is almost immediate. These are the cases of (mostly) Newton degenerate
singularities with reducible jets
(section \ref{SecReducibleJets}), i.e.$jet_p(f)$ defines a reducible hypersurface.
The simplest (nontrivial) example is the degenerate multiple point of order $p$ with hypersurfaces of a
form: $f(z_1,\dots,z_n)=\prod_i \Omega^{(p_i)}_i+higher~order~terms$.
Here $\Omega^{(p_i)}_i$ are non-degenerate mutually generic homogeneous forms of orders $p_i$ such that $\sum_ip_i=p$.

In this case, the equisingular family is defined by reducibility of the tensor of derivatives
(Appendix A). As reducibility is not invariant under topological transformations,
the equsingular family in this case does not coincide with the ND-topological stratum.
The enumeration goes in the same way as the enumeration of curves \cite{Ker}.
\paragraph{\hspace{-0.3cm}Computer calculations and efficiency.\hspace{-0.3cm}}
As we are working with polynomials of high degree in many variables, computer is used.
The calculations are essentially polynomial algebra: add/subtract polynomials, multiply,
open the brackets, eliminate variables, solve a big system linear equations etc.
Therefore a restriction arises from a computer's speed and memory.
We discuss some aspects of this step in $\S$ \ref{SecComputerCalcul}.

The restrictions are not severe for linear singularities, but are quite tough for non-linear ones.
In particular for $A_8$ (a "very non-linear" case) there was just not enough memory even in the case of curves.
We emphasize, however that this is a purely computer restriction.
\subsubsection{The simplest examples}
The case of ordinary multiple point is elementary because the defining conditions of the lifted stratum are globally
transversal. Usually the best we can hope for is to obtain the locally transversal conditions.

In this case, every time we degenerate, we should check for possible residual varieties and remove their
contributions if necessary. This technique is most simply demonstrated by the case of a node, naively defined in
affine coordinates.
\paragraph{Globally complete intersections: ordinary multiple point, $f=\sum z^{p+1}_i$.}\label{SecOrdinaryMultiplePoint}
We work with hypersurfaces $\{f(x)=0\}\subset\mP^n_x$ of degree $d$. The defining condition
here is: all the derivatives up to order $p$ should vanish. This can be written as $f|_x^{(p)}=0$
(tensor of derivatives of order $p$, in homogeneous coordinates, calculated at
the point $x$). The lifted variety in this case is:
\beq\label{node}
\tSi(x)=\{(f,x)|~f|_x^{(p)}=0\}\subset\mPN\times\mP^n_x
\eeq
(Recall, we always speak about topological closure of $\Si_\mD$.)
This variety is defined by ${p+n\choose{n}}$ transversal conditions. The transversality is proven in
general in $\S$ \ref{SecDefinCondLinearQuasiHomog}. For pedagogical reasons we check it here explicitly.

Note, that $\mP GL(n+1)$ acts freely and transitively on $\mP^n_x$,
therefore it is sufficient to check the transversality, at some particular point. For example
 fix $x=(1,0,\dots,0)\in\mP^n_x$.
Then the conditions of (\ref{node}) are just linear equations in the space $\mPN$ of all polynomials of the given
degree, so the transversality is equivalent to linear independence. And it is checked directly (note that $d$ is sufficiently high).

Thus the variety is a globally complete intersection and its cohomology class is just the product of
the classes of defining
hypersurfaces. Since all the hypersurfaces have the same class, we obtain\:
\beq\label{nodefull}
[\tSi(x)]=\left((d-p)X+F\right)^{n+p\choose{p}}
\eeq
Here $F,X$ are the generators of the cohomology ring of  $\mPN\times\mP^n_x$.

To obtain the cohomology class of $\Sigma$ we apply the Gysin homomorphism (as explained in
section \ref{SecMinimalLifting}).
From the expression in (\ref{nodefull}) one should extract the maximal non-vanishing power of $X$, i.e. $X^n$. The coefficient of this term is
the cohomology class of the needed stratum. This gives the degree:
\beq\label{DegMultPoint}
deg(\Sigma)={{n+p\choose{p}}\choose{n}}(d-p)^n
\eeq
\paragraph{Locally complete intersections: nodal hypersurfaces defined in affine coordinates.}
Let $x=(z_0\dots z_n)$ be the homogeneous coordinates in $\mP^n_x$.
Choose the affine part: $(z_0\neq0)\subset\mP^n_x$. In local coordinates a hypersurface has a node if the corresponding
function vanishes together with its derivatives. Thus we define:
\beq\label{nodeaff}
\widetilde\Xi(x)=\{(f,x)|\di_1f|_x=\dots=\di_nf|_x=0,~f|_x=0\}\subset\mPN\times\mP^n_x
\eeq
Over the affine part $\{z_0\ne0\}\subset\mP^n_x$ this variety coincides with the lifted stratum of nodal curves,
$\tSi_{A_1}(x)$.
 However at infinity ($z_0=0$) one can expect some additional pieces.
 Indeed, the Euler formula (\ref{Euler}) shows that the equations of (\ref{nodeaff}), when translated to
 the neighborhood of infinity, are:
\beq
z_0\di_0f|_x=0=\di_1f|_x=\dots=\di_nf|_x
\eeq
That is, the (projective closure of the) variety of (\ref{nodeaff}) is reducible: it is the union of
$\tSi_{A_1}(x)$ and some residual variety (at $z_0=0$), taken with multiplicity one
(since $z_0$ appears in the first degree).

In terms of co-homology classes:
\beq
[\widetilde\Xi(x)]=[\tSi_{A_1}(x)]+1[z_0=0,~\di_1f=\dots=\di_nf=0]
\eeq
So, to calculate the (co)homological class
$[\tSi_{A_1}(x)]$ one should subtract from $[\widetilde\Xi(x)]$ the (co)homological class of the variety defined by:
$\{z_0=0,~\di_1f=\dots=\di_nf=0\}$.
 Explicit calculation gives the degree of the discriminant (i.e. the result of
(\ref{DegMultPoint}) in the case $p=1$).
\subsubsection{Organization of material}
The main body of the paper gives the proof of the corollary \ref{ClaimMainOfTheMethod}.

In $\S$ 2 we recall some important definitions, fix the notations and introduce some auxiliary
notions used throughout the paper. We discuss the singularity types and the strata (topological and ND-topological,
section \ref{SecOnEquisingularity}) and clarify their relation. In $\S$ \ref{SecLinearSingularDefinitions}
we  introduce
 {\it linear singularities}. Then (in \ref{SecDefiningConditionsGeneral}) we formulate
{\it covariant defining conditions} for linear singularities.

In $\S$ \ref{SecLiftings} we define the {\it liftings}
of the strata and consider related questions. Then discuss the problem of non-transversality
(\ref{SecPossibleCyclesOfJump}) and obtain the characterization of points of non-transversal intersection.

In $\S$ \ref{SecIdeologyOfDegenerations} we prove the main theorems (\ref{TheoremEnumLinearSings} and \ref{TheoremEnumNonLinearSings}).
Essentially just collecting all the results from sections \ref{SecDefinitionsAuxilliaryResults}, \ref{SecMainIssues} and Appendix.

In $\S$ \ref{SecExamples} we demonstrate the algorithm by simple examples: the cusp $A_2$  and
the tacnode $A_3$. In $\S$ \ref{SecFurtherCalculations} we consider some higher singularities. We start from
the double points of a given co-rank, this corresponds to $A_2,D_4,P_8..$
Then, further degenerations are considered (e.g. $A_4,D_5,E_6$).

In most part of the paper we work with \gNnd singularities. In Appendix A we
consider a special subclass of Newton-degenerate singularities: with reducible
jets.

Appendix B is devoted to some results from intersection theory, intensively used throughout the paper  (multiplicity of
intersections, cohomology classes of some special varieties and restrictions of fibrations).

In Appendix C we give some explicit results (cohomology classes) and discuss the issues of computer
calculations and some consistency checks of the formulae.
\subsubsection{Acknowledgements}
I wish to thank my supervisor E. Shustin for his constant help and support during the work. The solution of the problem
would have been impossible without his ideas and great patience, when answering my technical questions
or during endless rewriting of the paper.

This work was initiated during the conference: "Singularities and Computer Algebra" on the occasion of
Gert-Martin Greuel's 60'th birthday (at Universit\"at Kaiserslautern).
I would like to thank the organizers for the invitation, L.Bodnarchuk and I.Burban
for the hospitality and numerous important conversations.

The work was inspired by R.Piene, she encouraged me and gave important advices.

To perform the explicit numerical calculations, I had to use computer programs. I would
like to thank B.Noyvert for numerous computational tricks, they tremendously simplified the process.

I would like to thank the referees for carefull reading of the text, their comments helped to improve it significantly.

Special thanks are to G.M.Greuel and A.Nemethi, who saved me from lots of errors in singularity of hypersurfaces.

The research was constantly supported by the Hermann-Minkowski Minerva Center for Geometry at Tel-Aviv
University and by Israel Science Foundation grant, no: 465/04.

\section{Some definitions and auxiliary results}\label{SecDefinitionsAuxilliaryResults}
\subsection{On variables and notations}\label{SectVar}
\subsubsection{On coordinates.\hspace{-0.3cm}} In this paper we deal with many varieties, embedded into various (products of)
projective spaces. Adopt the following convention.
If we denote a point in a projective space by the letter $x$, the corresponding projective space is denoted by $\mP^n_x$.
The points of projective spaces will be typically denoted by $x,y$ or $y_i$. For enumerative purposes we use
homogeneous coordinates: $x=(z_0:\dots:z_n)\in\mP^n_x$. When considering a particular singular hypersurface-germ we use
local coordinates, centered at the singular point, e.g. $(z_1,\dots,z_n)$ (assuming $z_0=1$).

When working with multi-projective space, the point
$(x_1,\dots,x_k)\in\mP^n_{x_1}\times\dots\times\mP^n_{x_k}$ is called generic if no subset of the points
$x_{i_1},\dots,x_{i_l}$ lies in a $(l-2)$-plane. The points $(x_1,\dots,x_k)$ will always be assumed mutually generic,
unless a restriction is explicitly specified.
By identifying $\mPN=\mP roj(V)$, we will often consider a point $x\in\mP^n_x$ as a vector of $(n+1)$ dimensional
 vector space $V$ (defined up to a scalar multiplication). For example, the above condition of genericity can be formulated as:
the vectors $x_1,\dots,x_k$ are linearly independent.

A hyperplane in $\mP^n$ is defined by a 1-form: $l\in(\mP^n_l)^*$. So, e.g. the incidence variety of
hyperplanes and their points is defined as:
\beq
\{(l,x)|l(x)=0\}\subset\mP^n_x\times(\mP^n_l)^*
\eeq
\subsubsection{On the monomial order.\hspace{-0.3cm}}\label{SecOnTheMonomialOrder} For the purpose of degeneration we should fix an
order on monomials $\bf{z}^{\bf I}$. So, we say that $\bf{z}^{\bf I}<\bf{z}^{\bf J}$ if $|I|<|J|$ (the total degrees).
For $|I|=|J|$ the order could be defined quite arbitrarily, we chose the lexicographic: $z_1>z_2>...z_n$.
\subsubsection{On symmetric forms.\hspace{-0.3cm}}We will often work with
symmetric $p-$forms $\Omega^{p}\!\in\! S^p V^*$
(here $S^p V$ is a symmetric power of $(n+1)$ dimensional vector space).
Thinking of a form as being a symmetric tensor with $p$ indices ($\Omega^{(p)}_{i_1,\dots,i_p}$), we often write
$\Omega^{(p)}(\underbrace{x,\dots,x}_{k})$ as a shorthand for the tensor multiplied $k$ times by a point
$x=(z_0,\dots,z_n)\in\mP^n$ (considered here as a vector in $V$):
\beq
\Omega^{(p)}(\underbrace{x,\dots,x}_{k}):=\sum_{0\le i_1,\dots,i_k\le n}\Omega^{(p)}_{i_1,\dots,i_p}z_{i_1}\dots z_{i_k}
\eeq
So, for example, the expression $\Omega^{(p)}(x)$ is a $(p-1)$ form. Unless stated otherwise, we assume the
symmetric form $\Omega^{(p)}$ to be generic (in particular non-degenerate, i.e. the corresponding hypersurface
$\{\Omega^{(p)}(\underbrace{x,\dots,x}_{p})=0\}\subset\mP_x^n$ is smooth).

Symmetric forms will typically occur as tensors of derivatives of order $p$: e.g. $f^{(p)}$ (here $f$ is a homogeneous
polynomial defining a hypersurface). Sometimes, to emphasize
the point at which the derivative is calculated we assign it. So, e.g.
$f|_x^{(p)}(\underbrace{y,\dots,y}_k)$ means: the tensor of derivatives of
$p$'th order, calculated at the point $x$, and contracted $k$ times with $y$.

Throughout the paper we tacitly assume the Euler identity for a homogeneous
polynomial of degree $d$
\beq\label{Euler}
\sum_{i=0}^nz_i\di_if=df
\eeq
and its consequences (e.g. $\sum_iz_i\di_i\di_jf=(d-1)\di_jf$). So, for example, the nodal point,
 defined by $f|_x^{(1)}=0$, can also be defined by $f|_x^{(p)}(\underbrace{x,\dots,x}_{p-1})=0$.
\subsubsection{On cohomology classes.\hspace{-0.3cm}}The generator of the cohomology
ring of $\mP^n_x$ are denoted by the corresponding upper case letter $X$, so
that $H^*(\mP^n_x)=\mZ[X]/(X^{n+1})$. Alternatively, $X=c_1({\cal O}_{\mP^n_x}(1))$. By the same letter we also
denote the hyperplane class in homology of $\mP^n_x$. Since it is always clear,
where we speak about coordinates and where about (co)homology classes, no confusion should arise.

To demonstrate this, consider the hypersurface:
\beq
\Si=\{(x,y,f)|~f(x,y)=0\}\subset\mP^n_{x}\times\mP^n_{y}\times\mPN
\eeq
Here $f$ is a polynomial of bi-degree $d_x,d_y$ in homogeneous coordinates $x=(z_0:\dots:z_n),y=(w_0:\dots:w_n)$,
the coefficients of $f$ are the homogeneous coordinates of the parameter space $\mPN$.
The cohomology class of this variety is:
\beq
[\Si]=d_xX+d_yY+F\in H^2(\mP^n_{x}\times\mP^n_{y}\times\mPN)
\eeq

The formulae for the cohomology classes of the lifted strata are polynomials
in the generators of the cohomology rings
of the products of projective spaces ($X$ for $\mP^n_{x}$, $F$ for $\mPN$ etc..) The polynomials depend
actually on some combinations of the generators. For example, $F$ always enters as $F+(d-k)X$ (for some $k\in\mN$,
which depends on the singularity type only). Correspondingly, we always use the (relative) class
\beq
Q:=(d-k)X+F
\eeq
(the value of $k$ will be specified or evident from the context).
\subsubsection{On the strata.\hspace{-0.3cm}}
We denote an ND-topological stratum by $\Sigma$ (it will be always clear from the context, which singularity
type is meant). The lifted stratum is denoted by $\tSi$. Usually there will be many liftings for one stratum,
to distinguish between them, we assign the auxiliary parameters. So, e.g. the stratum defined in
(\ref{ExampleLiftings}) is denoted by $\tSi(x)$.

We always work with the topological closures of the strata.
To simplify the formulae we
 write just $\Sigma$ (or $\tSi$) for the closure of a (lifted) equisingular stratum.
The lifted stratum  $\tSi$ is often considered as a fibration over the auxiliary space. For a cycle $C$ in
the auxiliary space the (scheme-theoretic) restriction of the fibration to the cycle is denoted by: $\tSi|_C$. For example,
 if $x$ is a point of the auxiliary space, then $\tSi|_x$ is the fibre over $x$.
\subsection{On the residual varieties and cycles of jump}\label{SecResidualVarieties}
We  try to represent a stratum as an explicit intersection of hypersurfaces. The intersection will
be usually non-transversal. The resulting variety, being reducible, will contain (except for the needed stratum) some
additional pieces. We call these pieces {\bf residual varieties}.

The intersection process occurs in the space: $\mPN\times Aux$. Here  $\mPN$
is the parameter space of hypersurfaces, while $Aux$ is the {\bf auxiliary space}, used to define
 the lifted stratum explicitly. It will be typically a multi-projective
space $Aux=\mP^{n_1}\times\dots\times\mP^{n_k}$ or a projective irreducible subvariety of it.
 We often consider the lifted stratum as a fibration over the auxiliary space.
The fibres correspond to hypersurfaces with some specified parameters of
the singularity (e.g. chosen singular point, tangent cone etc.) The fibration will be generically
locally trivial (in Zariski topology) and this local triviality induces the stratification on the auxiliary space.
\bed
Let
$R_1\subset\bar{R}_1=Aux$ be the maximal (open, dense) subvariety of $Aux$ over which the fibration $\tSi\ra Aux$
is locally trivial.
By induction, let
$R_k\subset\bar{R}_{k-1}\backslash R_{k-1}=Aux\backslash(\bigcup_{i=1}^{k-1}R_i)$ be the maximal subvariety such
that the (scheme) restriction $\tSi|_{R_k}\rightarrow R_k$ is a locally trivial fibration (over each connected component).
The set $\{R_i\}_i$ is
called: {\bf the stratification of the auxiliary space by local triviality}.
\eed
Let $\bar{R}_k\backslash R_k=\bigcup m_i C_i$ be the decomposition to a union of closed irreducible
subvarieties of $Aux$, not containing each other.
We call these irreducible subvarieties: {\bf the cycles of jump}.
The generic fibres of the projection $\tSi\ra Aux$ will be linear spaces, therefore the dimension of fibers jumps
over the cycles of jump (cf. $\S$ \ref{SecTypesSingul}).

Various cycles of jump can contain or intersect other cycles of jump (that appeared from $R_k$ with higher $k$).
It is useful to introduce grading on these cycles.
\bed\label{DefCycleOfJumpGrading}
The cycle of jump $C_i$ is assigned grading 1 if it does not contain other (proper) cycles of jump (equivalently
if the fibration $\tSi|_{C_i}\rightarrow C_i$ is locally trivial). A cycle of jump is assigned grading $k$
if it contains a (proper) cycle of grading $(k-1)$ and no cycles of higher grading.
\eed

\bex Quadratic forms of co-rank $r$
(this example is important for enumeration of singularities as $A_2$,$D_4$,$P_8\dots$).

Consider the variety of degenerate symmetric matrices (quadratic forms) of co-rank $r\geq2$ acting on $\mP^n$.
It is a classical determinantal variety \cite{Ful}. Its lifted version is
the incidence variety of degenerate quadratic forms with $r$ vectors of the kernel specified.
\beq\ber
\tSi(x_1,\dots,x_r):=
\overline{\Bigg\{(\!\!\!\underbrace{x_1,\dots,x_r}_{\ber\mbox{do not lie in an}\\(r-2)\mbox{-plane}\eer}\!\!\!,\Omega^{(2)})\Big|
\Omega^{(2)}(x_1)=\dots=\Omega^{(2)}(x_r)=0\Bigg\}}\subset\mP^n_{x_1}\times\dots\times\mP^n_{x_r}\times\mP_\Omega
\eer\eeq
(here $\mP_\Omega$ is the parameter space of quadratic forms, $dim(\mP_\Omega)={n+2\choose{2}}-1$,
$Aux=\mP^n_{x_1}\times\dots\times\mP^n_{x_r}$).

The projection: $(x_1,\dots,x_r,\Omega^{(2)})\stackrel{\pi}{\mapsto}(x_1,\dots,x_r)$ is generically locally trivial fibration.
The dimension of the generic fiber is ${n+2\choose{2}}-1-\frac{2n+3-r}{2}r$.
The {\it cycles of jump} here are all the diagonals: $\{(x_{i_1},\dots,x_{i_k})$ are linearly dependent$\}$.
The cycles of jump of minimal and maximal grades are:
\beq\ber
C_{min}=\{x_1=\dots=x_r\}\subset Aux~~~\rm{codim}_{Aux}(C_{min})=n(r-1)
\\C_{max}=\{(x_1\dots,x_r)\mbox{ lie in an }(r-2)\rm{-plane}\}\subset Aux~~~\rm{codim}_{Aux}(C_{max})=n+2-r
\eer\eeq
\eex
The first important question about cycles of jump is the jump in the dimension of fibers.
\bed
Let $c\in C,~x\in Aux$ be generic points of the cycle of jump and of the auxiliary space. The {jump of fiber dimension}
for the cycle $C$ is $\Delta\rm{dim}_C:=\rm{dim}(\tSi|_c)-\rm{dim}(\tSi|_x)$.
\eed
In the example above the jumps of dimension are:
\beq
\Delta\rm{dim}_{Cmin}:=\frac{(2n+2-r)(r-1)}{2}-1,~~~~~\Delta\rm{dim}_{Cmax}:=n+1-r
\eeq
The total variety, $\tSi$, will be always irreducible, in particular of pure dimension, therefore
we have immediate:
\bcor\label{ClaimJumpOfDimensionVsCodimension}
The jump of dimension over a cycle of jump is less than the codimension of the cycle of jump:
$\Delta\rm{dim}_{C}<codim_{Aux}(C)$
\ecor

The restriction of the fibration $\tSi\rightarrow Aux$ to the cycles of jump will be the source of
residual varieties, therefore we are interested in the cohomology classes of such restrictions: $\tSi|_C$.
This question is considered in $\S$ \ref{SecCohomologyClassOfRestrictionFibration}. By now
we need a simple technical result:
\bprop\label{ClaimCycleOfJump&Hypersurf}
Let $C$ be a cycle of jump and $\{C_i\}_i$ all the cycles of jump that are not contained in $C$. There exists
a hypersurface in the auxiliary space that contains $C$ and does not contain any of $C_i$ (though it can intersect them).
\eprop
\subsection{On singularities}\label{SecTypesSingul}
For completeness we recall some notions related to singularities of hypersurfaces \cite{AVGL,Dim,Shustin}.
\\
\parbox{13cm}
{
For a given (singular) hypersurface ($f=\sum a_{\bf I}{\bf z}^{\bf I}$) the Newton polytope is defined as
the convex hull of the support of $f$ in $\mZ^n$, namely: $conv({\bf I}\in\mZ^n|~a_{\bf I}\ne0)$.
We always take a generic representative for a given singularity type.
Therefore we can assume that the polynomial $f$ contains the monomials of the form $x^{d_i}_i,$ for $d_i\gg0$, and
the polytope intersects all the coordinate axes. Such a germ is called {\it commode} (or convenient).
In addition we assume (due to a high degree) that the hypersurface does not contain any line.
The upper part of the Newton diagram is defined as $\Gamma_+:=conv\Big(\bigcup ({\bf I}+\mR^n_+)|~a_{\bf I}\ne0\Big)$.
The Newton diagram is defined by $\Gamma_f:=\di\Gamma_+$.
By the above assumption it is compact and consists of a finite number of top dimensional faces.
}
\begin{picture}(0,0)(-60,0)
\put(0,0){\vector(0,1){60}}\put(0,0){\vector(-1,-1){50}}\put(0,0){\vector(1,0){80}}
\put(5,55){$\hat{e}_1$}   \put(-45,-55){$\hat{e}_2$}  \put(70,-10){$\hat{e}_3$}
\thicklines
\put(0,40){\line(-1,-2){40}} \put(-40,-40){\line(5,2){100}}   \put(0,40){\line(3,-2){60}}
\put(-2.5,38){$\bullet$}  \put(3,42){$p$} \put(-41,-41){$\bullet$} \put(-45,-35){$q$}
\put(58,-2.5){$\bullet$} \put(60,-10){$r$}
\put(5,-10){$\bullet$}
\put(8,-8){\line(-1,6){8}}  \put(7,-7){\line(-3,-2){47}}  \put(8,-8){\line(6,1){53}}
\put(-15,-40){\tiny The Newton diagram of $T_{pqr}$}    \put(30,35){$\Gamma_f$} \put(30,28){\vector(-1,-1){10}}
\put(-15,-50){\tiny $f=z^p_1+z^q_2+z^r_3+z_1z_2z_3$}
\end{picture}
\\
The restriction of a polynomial $f$ to its Newton diagram: $f|_{\Gamma_f}$ is called {\it the principal part}.

\subsubsection{On the singularity types and strata}\label{SecOnEquisingularity}
The embedded topological type and ND-topological type were defined in the introduction
(definition \ref{DefNDEquisingularStratum}).
\bed\cite[section 3.4]{Shustin}\label{DefEquiStratumNormalForm}
\li
The singular hypersurface $\{f=0\}$ is called Newton-non-degenerate if the restriction of
the polynomial $f$ to every face (of every dimension) of its  Newton diagram is non-degenerate
(i.e. the truncated
polynomial has no singular points in the torus $(\mC^*)^k$).
\li The singular hypersurface $f=0$ is \gNnd if it can be brought to a \Nnd form by a locally analytic transformation.
\eed
For example, while the hypersurface $(x-y^2)^2+y^5$
is not Newton-non-degenerate, it is certainly generalized Newton-non-degenerate.
\\
A {\it topological type} is called \Nnd if at least one of its representatives is \Nnd (or \gNnd).
Otherwise it is called Newton-degenerate.
The following natural question seems to be open (as always, we assume the topological strata to be algebraic and
irreducible):

Let $\{f=0\}$ be the {\it generic} representative of a \Nnd topological type. Is $f$ \gNnd?
\\\\
All the simple and the uni-modal singularity types are Newton-non-degenerate \cite[chapter 1]{AVGL}.
The simplest examples of not \gNnd hypersurfaces are  $W_{1,p\geq1}^\sharp$ with
Milnor number $\mu=15+p$
and $S^\sharp_{1,p\geq1}$ with Milnor number $\mu=14+p$. Even if a topological type is represented by not \gNnd
hypersurface
(e.g. as $W_{1,p\geq1}^\sharp$ before), it is not clear whether {\it every} representative of this type is not \gNnd
(this is the converse of the previous question).

A way to check this in particular cases was pointed to me by G.M.Greuel.
\bprop\label{ClaimNormalFamily}
Let $\{f_\al=0\}$ be a family of hypersurfaces (with one singular point) with the following properties:
\li They all have the same (local embedded) topological type $\mS=\mS_{f_\al}$.
\li The family covers all the moduli. Namely, for every value of moduli for the type $\mS$,
there is a representative $\{f_\al=0\}$ in this family with this value of moduli.
\li Every member of this family is {\it not} \gNnd.
\\Then every representative of the type $\mS$ is not \gNnd.
\eprop
The proof is immediate by observation that the whole stratum is a collection of equi-modular orbits under
the group of locally analytic transformations and our family intersects every orbit.

For example, from Arnol'd's classification it follows that the topological types
$W_{1,p\geq1}^\sharp$ and $S^\sharp_{1,p\geq1}$ are Newton-degenerate
(since the modality in both cases is 2 and the families of normal forms indeed cover all the moduli).
\\\\
\bed
\li The singular hypersurface  $\{f=0\}$ is called semi-quasi-homogeneous (SQH) if by a
locally analytic transformation it can be brought to a \Nnd form whose principal part is quasi-homogeneous.
\li The singularity type is called quasi-homogeneous if it has a SQH representative
\eed
In this case there is a strong result by \cite{Saito71}:
\bprop
Let $f$ be a quasi-homogeneous polynomial of degree $d$ and weights $w_1..w_n$, defining
a singular hypersurface, of topological singularity type $\mS_f$. Let another representative of this type
(algebraic hypersurface) be defined by a SQH polynomial $g$. Then $g$ is semi-quasi-homogeneous of the same degree $d$
and weights $w_1..w_n$ as $f$.
\eprop
~\\\\
In case of curves the embedded topological type and its strata possess all the good properties. For example,
the generic representative of the \Nnd type is \Nnd and can be brought to the given Newton diagram
by locally analytic transformations. Therefore the ND-topological type often coincides with the embedded topological type
(the same for the closures of the strata). For high enough degrees of curves the strata are irreducible and
smooth in their interior.

For hypersurfaces the situation is much more complicated. By choosing big degree $d$ of hypersurfaces, the non-emptiness
of the strata is assured.
But the strata of embedded topological type can behave quite badly. Just to mention, the algebraicity of the strata has
not yet been proven in general (though it is known for quite a broad class of types \cite{Var}).
Even assumed to be algebraic, the topological strata can be singular and reducible \cite{Luengo}
 (for large degrees of hypersurfaces).
The associated Newton diagram can be  non-constant along the equisingular stratum (in the sense that the transformation
needed to achieve it, is not locally analytic but a homeomorphism) \cite[chapter 1, example 2.14]{Dim}.
The constancy of multiplicity along the equisingular stratum has up to now been proved for
semi-quasi-homogeneous singularities only \cite{Greuel,GreuelPfister1} (see also \cite{EG} for recent results).
\\
\\
The equisingular strata we work with (the ND-topological strata) are
chosen especially to possess all the good geometric properties.
\bprop\label{ClaimNDtopStratumIsAlgVariety}
For a given Newton diagram $\mD$ the stratum $\Si_\mD$ is a (non-empty) irreducible algebraic variety
\eprop

Indeed, the family of hypersurfaces with a specific diagram
(the fiber over the diagram) is defined by linear equations in the parameter space $\mPN$.
And then the ND stratum is obtained by the action of algebraic group (locally analytic transformations) on this fiber.

From the irreducibility we get that every invariant, defined in an algebraic way (e.g. sectional Milnor numbers $\mu^*$)
is semi-continuous along the ND-topological strata.
\\
\\
By construction, to every ND-topological type $\mD$ an embedded topological type $\mS$ is associated,
with the inclusion of the (closures of the) strata: $\Si_\mD\subset\Si_\mS$. A natural question is:
when the two types of strata coincide?

\beR
The simplest example of non-coincidence (or non-uniqueness of ND-topological type for a given embedded topological)
 is just the
case of curves: $z^p_1+z_1x^p_2+z^d_2,~~d\ge p+2$.
\eeR
So, to get the equality of strata $\Si_\mS=\Si_\mD$ one should demand, at least, the minimality of Newton diagram.
In more details, introduce the partial order on the set of Newton diagrams with the same topological type: by inclusion.
Call a diagram {\it minimal} if it is
not bigger than some other diagram corresponding to the same topological type. In general, it is not
clear whether the minimal Newton diagram is unique (up to a permutation of axes).
Even if it is unique, it is still unclear whether the two strata coincide.

A constructive way to compare the strata is by codimension.
\bed
\li The codimension of a local embedded topological type $\mS$ is the codimension of
the topological stratum $\Si_\mS$ in the space of its semi-universal deformation. It equals: $\tau-\sharp(moduli)$
(here $\tau$ is the Tjurina number).
\li The codimension of an ND-topological type $D$, is the codimension of the stratum $\Si_\mD$ in
the parameter space of the hypersurfaces $\mPN$.
\eed
Restrict the topological type to be \Nnd and assume that the topological stratum is an irreducible algebraic variety.
From algebraicity and irreducibility we get:
\bcor
For a \Nnd type $\mS$, if the topological stratum $\Si_\mS$ is algebraic, irreducible
and the codimension of the topological type $\mS$ equals to that of the ND-topological type $\mD(\mS)$ then
the (closure of the) strata coincide: $\Si_\mS=\Si_\mD$.
So, in this case the generic representative of the topological type $\mS$ is \gNnd.
\ecor
Another observation is the following.
\bprop
Let $f_\al$ be the family as in the proposition \ref{ClaimNormalFamily}.
Assume also that they all have the same Newton diagram. Then the corresponding topological
and ND-topological strata coincide.
\eprop
Using this criterion we get the coincidence of the strata for all singularities with number of moduli $\le2$
(in this case the corresponding families are classified in \cite{AVGL}).

\subsubsection{On the vector spaces associated to the Newton diagram}\label{SecVectorSpacesOfNewtonDiagram}
The Newton diagram $\mD$ defines a stratification of the tangent space at the origin $T_0\mC^n$ as follows.
Let $\{f=0\}$ be the generic (\gNnd) hypersurface with diagram $\mD$,
let $l$ be a line through the origin, let the degree of their intersection be $k_l:=deg(k\cap \{f=0\})$.
(We assume the hypersurface to be generic, in particular it does not contain lines.)
The tangent space is stratified: $T_0\mC^n=\bigsqcup_k\cU_k$, according to the intersection
degree $\cU_k:=\{l\in T_0\mC^n|k_l=k\} $. Take the topological closures of $\cU_k$ and consider the irreducible components:
$\bar\cU_k=\cup_j \tV_{k,j}$. Call the collection of these components $\tcV$.
\bex For a SQH hypersurface, restrict to the principal (quasihomogeneous) part. Then the so obtained varieties $\tV_{kj}$
are just the vector spaces of a flag of $T_0\mC^n$. The flag may be non-complete, if some weight of (quasi-homogeneous)
variables coincide.
\eex
\bprop
The varieties $\tV_{kj}$ are vector subspaces of $\mC^n$ of the form$\{z_{i_1}\!=\!..\!=\!z_{i_j}\!=\!0\}$ (i.e. some coordinate planes).
They satisfy the property: if $\tV_{kj}\subsetneqq\tV_{k'j'}$ then $k>k'$.
\eprop
\bpr   Let $z_i=\al_it$ be the parametrization of a line through the origin. Consider its intersection
with the hypersurface $\{f=\sum a_{\bf I}{\bf z}^{\bf I}=0\}$, whose Newton diagram is $\mD$.
So, we study the function restricted to the line: $\sum a_{\bf I}({\bf \al t})^{\bf I}$.
If the line is contained in $\cU_k$ then $\sum_{I_1+..+I_n<k} a_{\bf I}({\bf \al t})^{\bf I}=0$
This equation must be satisfied for arbitrary (generic) coefficients of $f$ and for all small values of $t$.
Thus it causes the system of monomial equations of the form: ${\bf\al}^{\bf I}=0$. So, the irreducible components
of $\bar\cU_k$ are coordinate planes (which do not include one another).
\epr
\parbox{13.5cm}
{We need a more refined stratification. Every top dimensional face $\si$ of the Newton
diagram defines a flag of vector spaces as follows.
Let $k_i$ denote the non-zero coordinate of the intersection points of the $(n-1)$ plane $Span(\si)$
 with the coordinate axes (so $k_i$ are not necessarily integers).  Apply permutation
$\sigma\in S^n\subset GL(n)$ on the coordinate axes to arrange: $k_n\le\dots\le k_1$. Define vector spaces
 $\{\tilde{V}_i\}$ inductively:
\beq
\tilde{V}_n:=span(\hat{e}_1,\dots,\hat{e}_n),~~~\tilde{V}_{i-1}:=
\Big\{\ber \tilde{V}_i,~\rm{if}~k_{i-1}=k_i\\span(\hat{e}_1,\dots,\hat{e}_{i-1}),~\rm{if}~k_{i-1}>k_i\eer
\eeq
}
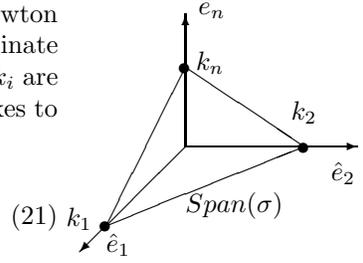
\begin{picture}(0,0)(-45,0)
\put(0,0){\vector(0,1){50}}\put(0,0){\vector(-1,-1){40}}\put(0,0){\vector(1,0){65}}
\put(5,50){$\hat{e}_n$}   \put(-30,-40){$\hat{e}_1$}  \put(55,-13){$\hat{e}_2$}
\put(0,-25){$Span(\si)$} 
\put(0,30){\line(-1,-2){30}} \put(-30,-30){\line(5,2){75}}   \put(0,30){\line(3,-2){45}}
\put(-3,27){$\bullet$}  \put(4,29){$k_n$}   \put(42,-3){$\bullet$}  \put(40,10){$k_2$}
\put(-33,-33){$\bullet$}  \put(-45,-30){$k_1$}
\end{picture}
\\
Take now the inverse permutation of axes ($\sigma^{-1}\in S^n\subset GL(n)$) that restores the initial values
of all $k_i$'s. Define: $V_i:=\sigma^{-1}(\tilde{V}_i)$.
\bed\label{DefFlagOfNewtonDiagrQuasiHomog}
The sequence of vector spaces: $\mC^n=V_n\supseteq\dots\supseteq V_1\supset\{0\}$ is called
{\bf the flag} of the face $\si$.
\eed
For a general (not SQH) singularity the Newton diagram consists of several top-dimensional faces,
each of them defines the corresponding flag. Now combine the flags together into the collection of
vector spaces $\{V_\al\}_{\al\in\cV}$. (The coinciding spaces are identified.)
We call $\cV$: {\it the collection of vector spaces associated to the Newton diagram}.
\bex \li For SQH hypersurface the collection is just a flag (as in the previous example).
\li The hypersurface $z_1^{p_1}+z_2^{p_2}+z_3^{p_3}+z_1z_2z_3+\sum_{i\ge4}z^2_i$ with
$\frac{1}{p_1}+\frac{1}{p_1}+\frac{1}{p_1}<1$, $p_1\lneqq p_2\lneqq p_3$ and $p_i>3$.
\\The top
dimensional faces are $Span(z_1^{p_1},z_2^{p_2},z_1z_2z_3,z^2_4...,z^2_n)$ and similarly for $(2,3)$ and $(1,3)$.
\\The flag of the top dimensional face is
\beq
0\subset Span(\he_2)\subset Span(\he_1,\he_2)\subset Span(\he_1,\he_2,\he_3)\subset\mC^n
\eeq
similarly for $(2,3)$ and $1,3$.
The collection of vector spaces is:
\beq
\{Span(\he_i)\}_{i=1,2,3},~~ \{Span(\he_i\he_j)\}_{\ber i\ne j\\i,j\le3\eer},~~
Span(\he_1\he_2\he_3),~~ \mC^n
\eeq
\eex
\beR
The latter stratification is finer than that by the degree of the intersection. As an example,
for the hypersurface $z_1^{2p}+z_2^{2p}+z_3^{2p}+z^{q-1}_1z^{q-1}_2z^{q-1}_3(z_1^3+z_2^3+z_3^3)+z_1^{p-1}z_2^{p-1}$
with $2p-2>3q$ and $q>1$ we have $\tcV=\bigcup_i\{z_i=0\}$ while
$\cV=\bigcup_i\{z_i=0\}\cup\{z_2=0=z_3\}\cup\{z_3=0=z_1\}$.
Therefore in the following we work with the collection $\cV$.
\eeR
As we consider the hypersurfaces of arbitrary dimensions, it is important to
check how the collection $\cV$ varies with $n$. More precisely, suppose $f_1$
(with $mult(f_1)>2$) is stably equivalent to $f_2$, i.e. $f_2(z_1..z_{n+k})=f_1(z_1..z_n)+\sum_{i>n}z_i^2$.
The relation between $\cV_{f_1}$ and $\cV_{f_2}$ is described by the simple proposition:
\bprop
\li Let $\si_1..\si_r$ be the top dimensional faces of $\Gamma_{f_1}$. Then the top dimensional faces
of $\Gamma_{f_2}$ are: $\Big\{conv(\si_i,\he_{n+1}.,\he_{n+r})\Big\}_{i=1..r}$.
\li Let $\{V_\al\}_{\al\in\cV_{f_1}}$ be the collection of vector spaces associated to $\Gamma_{f_1}$.
Then the collection $\cV_{f_2}$ is obtained as: $\{V_\al\}_{\al\in\cV_{f_1}}\cup
\Big\{Span(V_\al,\he_{n+1}.,\he_{n+r})\Big\}_{\al\in\cV_{f_1}}$.
\eprop
\beR
All the vector spaces above are defined in local coordinates. For enumerative purposes we need their
counterparts in homogeneous coordinates. For this, embed every $V\subset\mC^n$ into $\mC^{n+1}$, by
$\mC^n\ni(z_1,\dots,z_n)\stackrel{i}{\mapsto}(0,z_1,\dots,z_n)\in\mC^{n+1}$ and define:
\beq
\mC^{n+1}\supset\mV:=span(i(V),\hat{e}_0)
\eeq
The corresponding flag $\{\mV_i\}$ (or the collection of vector spaces) will be of key importance for writing down
the covariant defining equations of the strata.
\eeR
\subsubsection{Linear singularities}\label{SecLinearSingularDefinitions}
Fix a Newton diagram and the corresponding ND-topological type $\mD$. Consider a (hypersurface)
representative of this type.
In general, to bring it to the fixed Newton diagram a locally analytic transformation is applied. Split it into
steps are done (in local coordinates):
\\$\bullet$ Move the singular point to the origin. Rotate around the origin to fix the needed tangent cone.
\\$\bullet$ Make the (purely) quadratic transformation: $\vec{z}\rightarrow\vec{z}+\vec\Omega^{(2)}$
(here $\vec\Omega^{(2)}$ is an $n$-tuple of homogeneous quadratic forms), to remove some monomials.
\li Make the (purely) cubic transformation \dots
\dots

For some singularity classes the needed Newton diagram is achieved by just the first step (linear transformations).
As in the case of curves \cite[section 3]{Ker}, such singularities are much simpler for enumeration purposes.
As is shown later, these strata can be lifted to varieties, defined
by equations linear in function or its derivatives, and therefore easy to work with.
\bed\label{DefLinearSingularities}
\li For a given ND-topological type $\mD$, a hypersurface singularity (with this singularity type) is called linear,
if it can be brought to the needed Newton diagram by linear transformations only.
\li An ND-topological type/stratum is called linear if its generic representative is linear. Otherwise
the type/stratum is called non-linear.
\eed
The simplest linear singularity is the ordinary multiple point (here the needed diagram is achieved just
by translation of the singular point to the origin).
\\
\parbox{15cm}
{There is an easy characterization of linear singularities via their Newton diagrams:
\bel\label{ClaimCriterionLinearity}
The Newton-non-degenerate singularity type is linear iff all the angles between any face of the Newton diagram
 and the coordinate hyperplanes have the slope  $\frac{1}{2}\leq\mbox{tg}(\alpha)\leq2$.
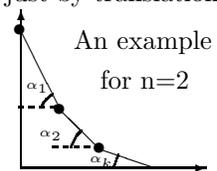
\eel
}
\begin{picture}(0,0)(-10,20)
\put(0,0){\vector(0,1){60}}\put(0,0){\vector(1,0){70}}
\put(-3,50){$\bullet$} \put(0,52){\line(1,-2){15}}
\multiput(13,23)(-4,0){4}{\line(-1,0){2}} \put(2,30){\tiny$\alpha_1$}  \qbezier(12,28)(8,25)(8,23)
\put(12,20){$\bullet$} \put(13,24){\line(1,-1){15}}
\multiput(30,8)(-4,0){5}{\line(-1,0){2}}  \put(7,12){\tiny$\alpha_2$} \qbezier(23,14)(21,13)(19,8)
\put(27,5){$\bullet$} \put(30,7){\line(3,-1){20}}
\put(26,2){\tiny$\alpha_k$} \qbezier(37,5)(36,3)(35,0)
\put(20,45){An example}\put(30,30){for n=2}
\end{picture}\\
\bpr$\Leftarrow$
Suppose all the slopes are bounded as above and a hypersurface germ has been brought to the given
Newton diagram by a chain
of locally analytic transformations. Start undoing these transformations to achieve the initial germ. Immediate
check shows that any nonlinear analytic transformation (without linear part) has no effect on the points under the Newton diagram.
Correspondingly the monomials of the initial polynomials that lie under the Newton diagram are restored by linear
transformation only. But it means that the germ could be brought to the Newton diagram by linear transformations only.
\\$\Rightarrow$
Suppose at least one of the angles (of the diagram) does not satisfy the condition $\frac{1}{2}\leq\mbox{tg}(\alpha)\leq2$. Then there exists
a quadratic shift of coordinates that changes the Newton diagram. Of course, such shift cannot be undone by linear transformations.
\epr
\beR\label{CorollarLinearSingMultPOrderDeterm2p}
As follows from the lemma, every Newton-non-degenerate singularity of multiplicity $p$ and order of
determinacy $k$, with $k\le 2p$, is linear.
\eeR
In case of plane curves there is only one angle for every segment of the diagram, correspondingly the condition on
the singularity to be linear is not too restrictive.
In the low modality cases the curve singularities brought to a Newton diagram by projective (linear)
transformation are (all the notations are from \cite{AVGL}):
\\$\bullet$ Simple singularities (no moduli): $A_{k\leq3},~~D_{k\leq6},~~E_{k\leq8}$
\\$\bullet$ Unimodal singularities: $X_9(=X_{1,0}),~J_{10}(=J_{2,0}),~Z_{k\leq13},~~W_{k\leq13}$
\\$\bullet$ Bimodal: $Z_{1,0},~W_{1,0},~W_{1,1},~W_{17},~W_{18}$

In case of hypersurfaces (of dimension$\geq3$) only a few singularities in each series can be linear.
The low modality cases are:
\\$\bullet$ Simple singularities (no moduli): $A_{k\leq3},~~D_{k\leq5},~~E_{6}$
\\$\bullet$ Unimodal singularities: $P_8,~X_9(=X_{1,0}),~Q_{10},~S_{11},~U_{12},~T_{p,q,r}~\{p,q,r\}\leq4$

For surfaces, there are some additional unimodal linear singularities:
$Q_{11},Q_{12}$,$S_{12},Q_{2,0}$,$S_{1,p\leq1}$, $U_{1,q\leq2}$, $Q_{17}$, $Q_{18}$, $S_{16}$,$S_{17}$

We emphasize, that there is infinity of linear singularities. More precisely: for every singularity type $\mD_1$,
there is a linear singularity $\mD_2$, which is adjacent to $\mD_1$ (e.g. one could take as $\mD_2$ an
ordinary multiple point of sufficiently high multiplicity).

Even if the singularity is non-linear, one can consider the collection of singular hypersurfaces that can be brought to
the given Newton diagram by {\it linear transformations} only (or projective
transformations in homogeneous coordinates). This defines a subvariety of the non-linear stratum:
{\it the linear substratum}. Our method, of course, enables to calculate the cohomology classes
of both the true strata and their linear sub-strata.
\subsection{Defining conditions of the singularities}\label{SecDefiningConditionsGeneral}
For Newton-non-degenerate singularities the defining conditions of a singular point are read from the Newton
diagram. Consider the points under the Newton diagram. The corresponding  monomials should be absent i.e. the
corresponding derivatives should vanish. To define the stratum, one has to write these conditions in a covariant
form. For linear singularities this can be done in especially simple way, since one should achieve
the covariance under the group of projective transformations (or linear in local coordiantes). Every condition has a form
$\{f^{(p)}_{i_1,\dots,i_p}=0\}$ and is transformed by $\mP GL(n+1)$ to $\{f^{(p)}(y_1,\dots,y_p)=0\}$, where $\{y_i\}$ are
some points of $\mP^n_{y_i}$ (regarded here as $(n+1)-$vectors).

First consider some simple examples.
\bex $\bullet$ An ordinary point of multiplicity $p$. Here in local coordinates we have: $f|_x=0$, $\di_1f|_x=\dots=\di_nf|_x=0$,
$\dots$, $\Big\{\di_{i_1}\dots\di_{i_{p-1}}f|_x=0\Big\}_{i_1\dots i_{p-1}}$

Passing to the homogeneous coordinates (and using Euler formula (\ref{Euler}) and its consequences) we get
the defining conditions in a covariant form: $f|_x^{(p-1)}=0$.
\\
\parbox{14cm}
{$\bullet$\label{ExamFlagOfKernelsHomogeneousForms} An ordinary point of corank $r$. Consider the singularity with
the normal form $\sum_{i=1}^{r}z^{p+1}_i+\sum^n_{i=r+1}z^p_i$. (For $p=2$ these singularities are $A_2,D_4,P_8..$).
The defining conditions are read directly from the Newton diagram:
\beq\ber
f|_x=0,~\di_1f|_x=\dots=\di_nf|_x=0,~\dots,~\Big\{\di_{i_1}\dots\di_{i_{p-1}}f|_x=0\Big\}_{i_1\dots i_{p-1}},
\\\Big\{\di_{i_1}\dots\di_{i_{p-1}}\di_jf|_x=0\Big\}_{\ber 1\le i_1\dots i_{p-1}\le n\\ j\le r\eer}
\eer\eeq}
\begin{picture}(0,0)(-50,0)
\put(0,0){\vector(0,1){50}}       \put(5,45){$\hat{e}_n$}  \put(-2,31){$\bullet$} \put(-10,30){$p$}
\put(0,0){\vector(1,2){20}}   \put(25,40){$\hat{e}_{n-1}$}   \put(11,24){$\bullet$} \put(20,27){$p$}
\put(0,0){\vector(2,1){45}}   \put(50,30){$\dots$}   \put(23,10){$\bullet$} \put(33,9){$p$}
\put(0,0){\vector(1,0){70}} \put(50,-10){$\hat{e}_{r+1}$}  \put(32,-2.2){$\bullet$}  \put(35,-10){$p$}
\put(0,0){\vector(1,-1){35}} \put(37,-43){$\hat{e}_{r}$}  \put(21,-26){$\bullet$}  \put(35,-27){$p+1$}
\put(0,-30){$\dots$}
\put(0,0){\vector(-1,-1){33}}      \put(-45,-38){$\hat{e}_1$}  \put(-26,-27){$\bullet$}  \put(-45,-16){$p+1$}
\put(-25,-23.5){\line(1,0){50}}   \put(-23,-23.5){\line(2,5){23}}  \put(0,34){\line(2,-1){15}}
\put(13,28){\line(5,-6){13}}    \put(24,15){\line(4,-5){12}}  \put(35,0){\line(-1,-2){12}}
\end{picture}
\\
To transform them to covariant form, introduce the flag of the Newton diagram (as defined
in \ref{DefFlagOfNewtonDiagrQuasiHomog}):
\beq
\mC^n=V_n=\dots=V_{r+1}\supset V_{r}=span(\hat{e}_1,\dots,\hat{e}_r)=\dots=V_{1}\supset\{0\}
\eeq
and the corresponding flag $\{\mV_i\}_i$ in $\mC^{n+1}$.
Then the conditions in homogeneous coordinates are:
\beq
f|_x^{(p-1)}=0,~~f|_x^{(p)}(y)=0,~~\forall y\in \mV_{r}\subset \mC^{n+1}
\eeq
\parbox{15.3cm}{$\bullet$ The $A_k$ point: $z^{k+1}_1+\sum^n_{i=2}z^2_i$. Here the flag is:
$\mC^n=V_n=\dots=V_2\supset V_1=span(\hat{e}_1)\supset\{0\}$. The conditions are:
\beq
f|_x=0,~\di_1f|_x=\dots=\di_nf|_x=0,~\Big\{\di^i_1\di_*f|_x=0\Big\}_{0<i<\frac{k+1}{2}+1},~
\Big\{\di^i_1f|_x=0\Big\}_{\frac{k+1}{2}+1\le i\le k}
\eeq}
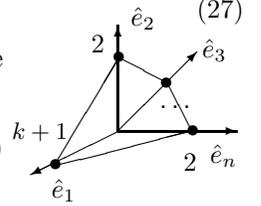
\begin{picture}(0,0)(-40,0)
\put(0,0){\vector(0,1){40}}       \put(5,40){$\hat{e}_2$}  \put(-2,26){$\bullet$} \put(-10,30){$2$}
\put(0,0){\vector(1,1){30}}   \put(32,28){$\hat{e}_3$}   \put(16,16){$\bullet$}  \put(16,9){$\dots$}

\put(0,0){\vector(1,0){45}} \put(35,-12){$\hat{e}_{n}$}  \put(26,-2.2){$\bullet$}  \put(25,-15){$2$}
\put(0,0){\vector(-2,-1){33}}      \put(-25,-25){$\hat{e}_1$}  \put(-26,-15){$\bullet$}  \put(-40,-3){\small$k+1$}

\put(-24,-13){\line(3,5){25}}   \put(-24,-13){\line(4,1){53}}
 \put(0,28){\line(2,-1){20}}  \put(28,0){\line(-1,2){10}}
\end{picture}
\\
$\bullet$ The $D_k$ point: $z^{k-1}_1+z^2_2z_1+\sum^n_{i=3}z^2_i$. We assume for simplicity that $k$ is even.
The flag is:  $\mC^n=V_n=\dots=V_3\supset V_{2}=span(\hat{e}_1,\hat{e}_2)\supset V_1=span(\hat{e}_1)\supset\{0\}$.
The conditions are:
\beq
f|_x=0,~\di_*f|_x=0,~~
\Big\{\di^i_2\di^{l-i-1}_1\di_*f|_x=0\Big\}_{\ber 0\le i<\frac{k+1-2l}{k-4}\\l<\frac{k+1}{2}\eer},~
\Big\{\di^i_2\di^{l-i}_1f|_x=0\Big\}_{\ber i<\frac{2(k-1-l)}{k-4}\\l<k-1\eer}
\eeq
\eex
We describe now the general procedure of formulating the covariant conditions for linear singularities.
Recall that in $\S$ \ref{SecVectorSpacesOfNewtonDiagram} the collection of vector spaces associated to a given
Newton diagram was defined. Start from the SQH case, where the collection is just a flag in $\mC^n$
(definition \ref{DefFlagOfNewtonDiagrQuasiHomog}).

\subsubsection{The case of semi-quasi-homogeneous linear singularities}\label{SecDefinCondLinearQuasiHomog}
\parbox{13.5cm}{
Let $\hat{e}_1,\dots,\hat{e}_n$ be the coordinate axes of the lattice (by the same letters we also denote the
corresponding unit vectors).
Consider the points of intersection of the hyperplane $\Gamma_f$
with the coordinates axes: $\{k_i\}_i$. By
renumbering the axes we can assume: $k_1\ge k_2\ge\dots\ge k_n$.

To write the conditions of the Newton diagram explicitly, consider the points of the
lattice lying under the Newton diagram. As in the examples above, at each step we consider
points corresponding to partial derivatives of a given order.
}
\begin{picture}(0,0)(-60,0)
\put(0,0){\vector(0,1){60}}\put(0,0){\vector(-1,-1){45}}\put(0,0){\vector(1,0){70}}
\put(5,50){$\hat{e}_n$}   \put(-35,-45){$\hat{e}_1$}  \put(60,-12){$\hat{e}_{n-1}$}
\put(-15,40){$\Gamma_f$} \put(-15,40){\vector(1,-2){10}}
\put(0,40){\line(-1,-2){40}} \put(-40,-40){\line(5,2){100}}   \put(0,40){\line(3,-2){60}}
\put(-20,-20){\line(3,1){60}}  \put(-20,-20){\line(1,5){6.5}}   \put(40,0){\line(-1,2){10}}
 \put(-13.5,13){\line(6,1){43}}  \put(35,20){$\Delta_r$} \put(35,20){\vector(-1,-1){10}}
\put(-3,36){$\bullet$}  \put(12,36){$k_n$}   \put(56,-3){$\bullet$}  \put(56,5){$k_{n-1}$}
\put(-42,-42){$\bullet$}  \put(-50,-35){$k_1$}
\end{picture}
\\
Namely, for every $0\le r\le k_1$ define an $(n-1)$ dimensional simplex:
\beq
\Delta_r:=\{x=(m_1,\dots,m_n)|~m_i\ge0,~\sum^n_{i=1} m_i=r,~~x\mbox{ lies strictly below~}\Gamma_f\}
\eeq
Every integral point of $\Delta_r$ corresponds to a vanishing derivative (in local coordinates):
$\di^{m_1}_1\dots\di^{m_n}_nf\equiv f^{(r)}_{\underbrace{\tinyM1\dots1}_{m_1},\dots,\underbrace{\tinyM n\dots n}_{m_n}}=0$.

We need a more precise relation between the axes of the diagram and the vector spaces
of the flag of the diagram.
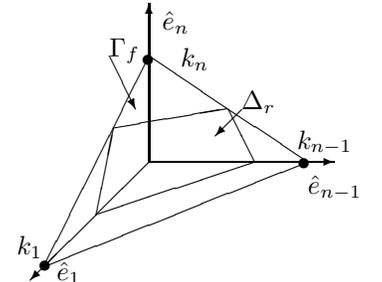
\bed
Let $\{V_i\}_i$ be the flag of the Newton diagram. For each axis $\hat{e}_i$ define the {\bf associated
vector space} as $V(\hat{e_i}):=V_j$, such that $\hat{e_i}\in V_j$ and $\hat{e_i}\notin V_{j-1}$ .
\eed
Recall also that each vector space $V_i\subset\mC^n$ has its (homogeneous) version $\mV_i\subset\mC^{n+1}$.

The transition from local conditions (that arise from the given Newton diagram) to the
conditions covariant under $PGL(n+1)$ is done by the following
\bel\label{ClaimDefiningConditionsQasHomSing}
For a SQH singularity the local set of conditions
$\Big\{\di^{m_1}_1\dots\di^{m_n}_nf=0\Big\}$ for $\Big\{(m_1\dots m_n)\in\Delta_r\Big\}$,
corresponding to the points under the Newton diagram, can be written in a covariant way as:
\beq
\forall y_i\in \mV(\hat{e}_i)~~~\Big\{\forall (m_1,\dots,m_n)\in \Delta_r\Big\}^{k_1}_{r=0}~~~~
f^{(r)}(\underbrace{y_1\dots y_1}_{m_1},\dots,\underbrace{y_n\dots y_n}_{m_n})=0
\eeq
\eel
\bpr $\Rightarrow$ We should check that the covariant conditions imply the local
conditions. This is immediate (just take: $y_i=\hat{e}_i$).
\\$\Leftarrow$ Without loss of generality, can assume $k_1\ge k_2\ge\dots\ge k_n$. Introduce the
lexicographic order on the points of $\Delta_r$ (for a fixed $r$):
\beq
(m_1,\dots,m_n)<(\tilde{m}_1,\dots,\tilde{m}_n)~~\rm{if}~\Big\{\ber m_i=\tilde{m}_i:~\rm{for}~i<j\\m_j<\tilde{m}_j\eer
\eeq
We have: if the point $(\tilde{m}_1..\tilde{m}_n)$ lies under the Newton diagram and
$(m_1..m_n)\leq(\tilde{m}_1..\tilde{m}_n)$, then the point $(m_1..m_n)$ also lies under
the Newton diagram (because: $\sum \frac{m_i}{k_i}\leq\sum \frac{\tilde{m}_i}{k_i}<1$).

Now, expanding: $y_i=\sum \alpha_i\hat{e_i}$ and substituting into
$f^{(r)}(\underbrace{y_1..y_1}_{m_1}..\underbrace{y_n..y_n}_{m_n})$ we get the sum of terms,
corresponding to the points $(\tilde{m}_1..\tilde{m}_n)$, satisfying:
$(m_1..m_n)\ge(\tilde{m}_1..\tilde{m}_n)$. Therefore, if $(m_1..m_n)\in\Delta_r$ then
$f^{(r)}(\underbrace{y_1..y_1}_{m_1}..\underbrace{y_n..y_n}_{m_n})=0$.
\epr
\bex\label{ExamFlagOfKernelsTacnode} Below we write the conditions in several cases, considered in the
example \ref{ExamFlagOfKernelsHomogeneousForms}.
\li Linear singularities of $A_k$ series are $A_{k\le 3}$. In the $A_2$-case, the covariant conditions are:
$f|_x^{(1)}=0$,  $\forall y\in \mV_1: f|_x^{(2)}(y)=0$. In the tacnodal case ($A_3$) the additional
condition is $f|_x^{(3)}(y,y,y)=0$.
\li Linear singularities of $D_k$ series are $D_{k\le 5}$.The $D_4$-case
was considered in example \ref{ExamFlagOfKernelsHomogeneousForms}. Here we consider $D_5$.
The covariant conditions are: $f|_x^{(1)}=0$, $\forall~y_1\in\mV_1,~y_2\in\mV_2:$ $f|_x^{(2)}(y_1)=0=f|_x^{(2)}(y_2)$,
$f|_x^{(3)}(y_1,y_1,y_1)=0=f|_x^{(3)}(y_1,y_1,y_2)$.
\li The singularities of ND-topological type with the representative:
$\sum^{r_1}_{i=1}z^{p+2}_i+\sum^{r_2}_{i=r_1+1}z^{p+1}_i+\sum^{n}_{i=r_2+1}z^p_i$. For $p=2$ this class
contains: $A_{k\le 3},$ $D_4,$ $E_6~(r_1+1=r_2=2),$ $P_8~(r_1=0,~r_2=3)$ $X_9~(r_1=r_2=2),$ $U_{12}~(r_1+2=r_2=3)$ etc.
\\
The flag is: $\mC^{n}=V_n=\dots=V_{r_2+1}\supset V_{r_2}=span(\hat{e}_{1},\dots,\hat{e}_{r_2})=\dots=
V_{r_1+1}\supset V_{r_1}=span(\hat{e}_{1},\dots,\hat{e}_{r_1})=\dots=V_1\supset\{0\}$.
\\
The point $(m_1,\dots,m_n)\in\Delta_r$ lies under the Newton diagram provided:
\beq
\sum_i m_i=r,~~~~~\sum^{r_1}_{i=1}\frac{m_i}{p+2}+\sum^{r_2}_{i=r_1+1}\frac{m_i}{p+1}+
\sum^{n}_{i=r_2+1}\frac{m_i}{p}<1
\eeq
The covariant equations are:
\beq\ber
f|_x^{(p-1)}=0,~\forall~y_1\in V_{r_1},~\forall~y_2\in V_{r_2}:~f|_x^{(p)}(y_2)=0,~
f^{(p+1)}|_x(\underbrace{y_1,\dots,y_1}_{k},\underbrace{y_2,\dots,y_2}_{l})=0,~\rm{for}~\frac{2k}{p+2}+\frac{l}{p+1}>1
\eer\eeq
\li The singularities of ND-topological type with the representative
$\sum^{r_1}_{i=1}z^{4}_i+z_{r_1}z^2_{r_1+1}+\sum^{r_2}_{i=r_1+2}z^{3}_i+\sum^{n}_{i=r_2+1}z^2_i$.
This class contains: $D_4,$ $P_8,$ $Q_{10},$ $V_{1,0}$ etc.
\\
The flag is: $\mC^{n}=V_n=\dots=V_{r_2+1}\supset V_{r_2}=span(\hat{e}_{1},\dots,\hat{e}_{r_2})=\dots=
V_{r_1+2}\supset V_{r_1+1}=span(\hat{e}_{1},\dots,\hat{e}_{r_1+1})\supset V_{r_1}=span(\hat{e}_{1},\dots,\hat{e}_{r_1})
=\dots=V_1\supset\{0\}$.
\\
The point $(m_1,\dots,m_n)\in\Delta_r$ lies under the Newton diagram provided:
\beq
\sum_i m_i=r,~~~~~\sum^{r_1}_{i=1}\frac{m_i}{4}+\frac{3m_{r_1+1}}{8}+\sum^{r_2}_{i=r_1+2}\frac{m_i}{3}+
\sum^{n}_{i=r_2+1}\frac{m_i}{2}<1
\eeq
The covariant equations are:
\beq\ber
f|_x^{(1)}=0,~\forall~y_1\in V_{r_1},~\forall~y_2\in V_{r_2},~\forall~y_3\in V_{r_1+1}:~f|_x^{(2)}(y_2)=0,~
f|_x^{(3)}(y_1,y_2,y_3)=0
\eer\eeq
\eex
\subsubsection{The general linear case}\label{SecDefinCondLinearNonQuasHom}
In the linear non-SQH case we have the collection of vector spaces for the top-dimensional faces of the Newton diagram.
So, perform the procedure for each top dimensional face separately (i.e. for the hyperplane that contains it).
Thus the lemma \ref{ClaimDefiningConditionsQasHomSing} is translated verbatim to a collection of lemmas, for each face.
Now combine all the conditions (the coinciding conditions should be identified).

So, this generalizes the method of obtaining defining conditions to the case of non-quasi-homogeneous singularities.

\bex
The singularity of the type $T_{pqr}$ with the normal form: $z^p_1+z^q_2+z^k_3+z_1z_2z_3+\sum_{i=4}^nz^2_i$
$(p\ge q\ge k)$. It is linear for $p,q,r\le5$. Here the Newton diagram consists of the three planes:
\beq\ber
\tinyT\tinyM L_{pq}:=\Big\{\!\!\frac{m_1}{p}+\frac{m_2}{q}+m_3(1-\frac{1}{p}-\frac{1}{q})+\sum_{i\ge4}\frac{m_i}{2}=1\!\!\Big\},~
\tinyM L_{pk}:=\Big\{\!\!\frac{m_1}{p}+m_2(1-\frac{1}{p}-\frac{1}{k})+\frac{m_3}{k}+\sum_{i\ge4}\frac{m_i}{2}=1\!\!\Big\},~
\tinyM L_{qk}:=\Big\{\!\!m_1(1-\frac{1}{k}-\frac{1}{q})+\frac{m_2}{q}+\frac{m_3}{k}+\sum_{i\ge4}\frac{m_i}{2}=1\!\!\Big\}
\eer\eeq
So, one considers the three flags: $(V_{i,pq},V_{i,qk},V_{i,pk})$ and the three sets of polytopes:
$(\Delta_{r,pq},\Delta_{r,qk},\Delta_{r,pk})$. We get a set of local conditions, which is transformed to the set of
covariant conditions, by the prescription of lemma \ref{ClaimDefiningConditionsQasHomSing}. We omit the calculations.
\eex
\subsubsection{On the transversality of covariant conditions}
We would like to address here the issue of transversality. All the initial conditions (corresponding
to the points under the Newton diagram) are just linear equations on
the (derivatives of the) function $f$, the transversality is equivalent to the linear independence. And this is
obvious for conditions that are read from the Newton diagram.
The conditions are made covariant by introducing vectors from the collection of vector spaces of the diagram.
As the linear independence is preserved
under (invertible) linear transformations we obtain that the covariant conditions remain transversal
as far as these vectors are mutually generic. More precisely:
\bprop\label{ClaimTransversalityUnderNewtonDiagram}
The set of conditions $\{f^{(p_j)}(y_{i_{1,j}},\dots,y_{i_{p_j,j}})=0\}_{i,j}$ (that correspond to the points
under the Newton diagram) is transversal if by $PGL(n+1)$ the points $\{y_{i,j}\}_i$ can be brought
to the coordinate axes.
\eprop
In particular one demands that no subset of $k$ variables lies in a $(k-2)$ hyperplane (for any $k$).

\subsubsection{The non-linear singularities and defining conditions}\label{SecDefinCondNonLinearSing}
Singularity types/strata for which a Newton diagram cannot be achieved by projective (linear) transformations are non-linear.
Most types/strata are nonlinear.
\bex\label{ExampleDefiningEqsA4}
Consider a hypersurface with the $A_4$ point. Try to bring it to the Newton diagram of $A_4$ (i.e.
that of $z^4_1+z^2_2+..z^2+n$).
The best we can do by projective
transformations is to bring it to a form of $A_3$:
\beq
f=\sum_{i=2}^n\alpha_iz_i^2+z_1^2\sum_{i=2}^n\beta_iz_i+\gamma z^4_1+\dots
\eeq
To achieve the Newton diagram of $A_4$ we must do the non-linear shift: $z_i\rightarrow z_i+\delta_iz^2_1$
(to kill the monomials $z^4_2,z^2_1z_i$). Elimination of the parameters of the transformation gives the
non-linear equation: $\gamma=\sum\frac{\beta^2_i}{4\alpha_i}$.
\eex

In general, to obtain the locally defining equations of a non-linear singularity one considers the
locally analytic transformation of coordinates:
\beq
z_i\rightarrow z_i+\sum_{k_i=1}^\infty a_{i,k}z_{k_i},~~~\rm{for}~~i=1\dots n
\eeq
Then, demanding that some derivatives (corresponding to the points under the Newton diagram) vanish, one
eliminates the parameters of the transformation $\{a_{ik}\}_{i,k}$ to obtain the system of (non-linear)
equations.
This forms the locally defining ideal of the stratum
The ideal can be quite complicated (since the stratum can be not a locally complete intersection).
It is important to trace the whole ideal, i.e. all of its generators. Of course the whole procedure is done
by computer and is time consuming.

Note that, as we enumerate the non-linear strata by the {\it simplifying degenerations}
(defined in \ref{SecSimplifyingDegenerations}), we do not need to convert the defining conditions to
a covariant form.

Among the generators of the ideal there are (a finite number of) non-linear expressions in the coefficients of $f$.
The goal of degeneration is to turn these non-linear equations into monomial equations. In this case
the ideal would correspond to a collection of linear strata (with some multiplicities).

For enumeration of a non-linear singularity we will be interested in the set of linear singularities to which the
given singularity is adjacent.
\bed
The linear singularity type $L$ is called assigned to a given singularity type $N$ if $N$ is adjacent to $L$
(i.e. $\Si_N\subset\bar{\Si_L}$) and the adjacency is minimal. Namely there is no other linear type $L'$ such that
$\Si_N\subset\bar{\Si_L'}\subsetneqq\bar{\Si_L}$.
\eed
The assigned type is non-unique in general.
As singularities always appear in series, which start from linear singularities, we have a natural way to
fix an assigned linear singularity. In the simplest cases the assigned linear singularities are:
 $A_3$ for $A_{k\ge4}$, $D_5$ for $D_{k\ge6}$, $E_6$ for $E_{k>6}$ etc.
\subsubsection{On the invertibility of degenerations}\label{SecInvertibilityofDegeneration}
At each step of the degeneration procedure we get an equation in the cohomology ring of $Aux\times\mPN$:
\beq
[\tSi_1][\mbox{degenerating divisor}]=[\tSi_2]
\eeq
An important issue is to check that the degeneration is "invertible", i.e. this equation
 fixes the cohomology class of the original stratum uniquely. The generators of the cohomology ring are nilpotent
(e.g. $X^{n+1}=0=F^{D+1}$). Thus the solution for $[\tSi_1]$ is unique, provided
$dim(\Si_1)+1\le D$ and the class of degenerating divisor depends in essential
way on $F$.
The first condition is always satisfied, while the second means that the degeneration must involve (in essential way)
the function  $f$ (or its derivatives). In particular, conditions involving the parameters of the auxiliary space only
(e.g. coincidence of points) are non-invertible and will not be used for the degenerations.
\section{Enumeration}
\subsection{Main issues of the method}\label{SecMainIssues}
\subsubsection{Liftings (desingularizations)}\label{SecLiftings}
We lift the strata to a bigger ambient space to define them by explicit equations.

The simplest example (minimal lifting),
consisting of pairs (the function, the singular point), was considered in \ref{SecMinimalLifting}.
This lifting is sufficient for ordinary multiple points only. In general one should lift further
to the space  $Aux\times\mP_f^D$. We will
consider the objects:
(the function, the singular point, hyperplanes in the tangent cone, some special linear subspaces of the hyperplanes).

In $\S$ \ref{SecVectorSpacesOfNewtonDiagram} we constructed for a given Newton diagram $\mD$ the
collection of vector spaces (in homogeneous coordinates) $\{\mV_\al\}\in\cV(\mD)$. Taking their projectivization
we arrive at a collection $\{\mP_{y_\al}^{n_\al}\}_{\al\in\cV}$ (with the conditions that some
are subspaces of others). This defines the lifting:
\bed\label{DefLifting} For linear ND-topological type the lifting is defined by
\beq
\tSi(x,\{y_\al\}_i):=\overline{\Bigg\{\!\!\!
\ber (x,\{y_\al\},f)\\\mbox{generic}\eer\!\!{\Bigg|}
\!\!\ber\Big\{y_\al\in\mP_{y_\al}^{n_\al}\Big\}_\al\in{\cV}(\mD),~~f
  \mbox{ has the prescribed}\\\mbox{ singularity with }
\mD~\mbox{as its Newton diagram}\eer\Bigg\}}\subset
\mP^n_x\times\prod_\al\mP_{y_\al}^{n_\al}\times\mPN
\eeq\eed
\bex In the SQH case the collection of vector spaces is just the flag therefore
\beq
\tSi(x,\{y_i\}_i):=\overline{\Bigg\{\!\!\!
\ber (x,\{y_i\}_i,f)\\\mbox{generic}\eer\!\!{\Bigg|}
\!\!\ber\Big\{\mV_i=span(y_1,\dots,y_{dim(\mV_i)})\Big\}_i\!\!\!\!\\\ f \mbox{ has the prescribed singularity with}\\
\mV_n\supset\dots\supset\mV_{1}~\mbox{as the flag of the Newton diagram}\eer\Bigg\}}\subset
\mP^n_x\times\prod_i\mP^n_{y_i}\times\mPN
\eeq
\eex
We emphasize that the tuples $(x,\{y_i\}_i)$ are always taken to be generic (in particular, no $k$ of the points
span a $(k-2)$-plane). For the sake of exposition we will often omit this, but it is always meant.
\beR
As follows from lemma \ref{ClaimDefiningConditionsQasHomSing} the so defined stratum $\tSi$ is indeed a lifting.
That is: every point of $\tSi$ projects to a point of $\Si$.
This definition reduces the enumerative problem for linear singularities to the intersection theory.
The intersection theory step has still many complications, related questions are considered in Appendix B.
\eeR
\bex\label{ExampleNumerousLiftings} For many cases the defining covariant conditions were given in example \ref{ExamFlagOfKernelsTacnode}.
Therefore, we immediately have:
\\$\bullet$ For multiple point of co-rank $r$ (with normal form $f=\sum_{i=r+1}^n x_i^p+higher~order~terms$).
\beq\ber\label{CorankLiftings}
\tSi(x,y_1,\dots,y_r)=\overline{\Big\{\ber(x,y_1,\dots,y_r)\\~~generic\eer|
f|_x^{(p-1)}=0=\Big(f|_x^{(p)}(y_i)\Big)_{i=1}^r\Big\}}\subset\mP^n_x\times\prod_{i=1}^r\mP^n_{y_i}\times\mPN
\eer\eeq
\\$\bullet$ For singularities with the normal form:
$\sum^{r_1}_{i=1}z^{p+2}_i+\sum^{r_2}_{i=r_1+1}z^{p+1}_i+\sum^{n}_{i=r_2+1}z^p_i$:
\beq
\tSi(x,y_1,\dots,y_{r_2})=\overline{\Big\{(x,y_1,\dots,y_{r_2)}\Big|\ber
f|_x^{(p-1)}=0~f|_x^{(p)}(y_1)=0~\rm{for}~i\le r_2\\
f^{(p+1)}|_x(y_{i_1},\dots,y_{i_k},y_{j_1},\dots,y_{j_l})=0\\
~\rm{for}~i_1..i_k\le r_1,~j_1..j_l\le r_2~~\frac{2k}{p+2}+\frac{l}{p+1}>1\eer\Big\}}
\subset\mP^n_x\times\prod_{i=1}^{r_2}\mP^n_{y_i}\times\mPN
\eeq
$\bullet$ For singularities with the normal form:
$\sum^{r_1}_{i=1}z^{4}_i+z_{r_1}z^2_{r_1+1}+\sum^{r_2}_{i=r_1+2}z^{3}_i+\sum^{n}_{i=r_2+1}z^2_i$:
\beq
\tSi(x,y_1,\dots,y_{r_2})=\overline{\Big\{(x,y_1,\dots,y_{r_2)}\Big|\ber
f|_x^{(1)}=0~f|_x^{(2)}(y_i)=0~\rm{for}~i\le r_2\\
f|_x^{(3)}(y_i,y_j,y_k)=0~\rm{for}~i\le r_1,~j\le r_1+1,~k\le r_2\eer\Big\}}
\subset\mP^n_x\times\prod_{i=1}^{r_2}\mP^n_{y_i}\times\mPN
\eeq
\eex
\paragraph{\hspace{-0.3cm}Lifted varieties as fibrations over the auxiliary space.}\hspace{-0.3cm}~
\\
Note that for a given singular germ $(f=0)$ only the flag $\mV_n\supset\dots\supset \mV_{1}$ and
the assigned polytopes are defined uniquely. The points $\{y_i\}_i$ of the auxiliary space
(in definition \ref{DefLifting}) can vary freely as far as the flag and the polytopes
are preserved. So we obtain an important result:
\bprop
The projection: $\tSi\stackrel{\pi}{\rightarrow}\Sigma$ is a fibration with generic
fiber the multi-projective space.
\eprop
In the example of multiple point of co-rank $r$ the generic fiber is $\mP^r_{y_1}\times\dots\times\mP^r_{y_r}$.

Note, that the nonzero dimensionality of the fibers already restricts the possible cohomology class of the lifted variety.
In fact, let $(X,Y_1,\dots,Y_r,F)$ be the generators of the cohomology rings of $\mP^n_x,\mP^n_{y_i},\mPN$
(as defined in $\S$ \ref{SectVar}).
The cohomology class of $\tSi$ is a polynomial in $(X,Y_1,\dots,Y_r,F)$.
\bel\label{ClaimFirstConsistencyCondition} {\bf (The first consistency condition).}\\
\parbox{12cm}
{Let the fibration with generic fiber $\mP^r_{y_k}\subset\mP^n_{y_k}$ be given as in the diagram.
Then the variable $Y_k$ that appears in monomials of the polynomial $[\tSi(x,y_1,\dots,y_k)]$
has powers not bigger than $(n-r)$.}~~~~~~
$\ber
\tSi(x,y_1,\dots,y_k)\subset Aux\times\mP^n_{y_k}\times\mP^D
\\~~~~~~~~~~ \downarrow~~~~~~~~~~~~~~~~\downarrow\\
\tSi(x,y_1,\dots,y_{k-1})\subset Aux\times\mP^D
\eer$
\eel
\parbox{15cm}
{\bpr We give the proof in the general case. Consider projective fibration $E\ra B$ with the generic fiber $\mP^r_y$.
Assume that the generic fiber is linearly embedded into $\mP^n_y$.
Write the cohomology class $[E]\in H^*(A\times \mP^n_y)$ as:
}~~~
$\ber
E\hookrightarrow A\times \mP^n_y\\\downarrow~~~~~\downarrow\\B\hookrightarrow A
\eer$
\beq
[E]=\sum_i Y^iQ_{dim(A)+n-dim(E)-i}
\eeq
here $Y$ is a generator of $H^*(\mP^n_y)$, while $Q_{dim(A)+n-dim(E)-i}\in H^{2(dim(A)+n-dim(E)-i)}(A)$. We want to
show that terms appearing in the above sum have: $i\le n-r$. To see this, multiply $[E]$ by
$Y^{n-i}\tilde{Q}_{dim(E)-n+i}$ (for some arbitrary $\tilde{Q}$). By the duality between homology and cohomology this product corresponds
to the intersection of $E$ with (generic) cycle of the form:
$\mL^{(i)}_{\mP^n}\times C^{dim(A)-dim(E)+n-i}_A$ (here $\mL^{(i)}_{\mP^n}$ is a linear $i-$dimensional
subspace of $\mP^n$, $C-$a cycle in $A$).
Then (by the dimensional consideration in $A$) the intersection is empty unless:
\beq
dim(A)-dim(E)+n-i+dim(B)\ge dim(A)
\eeq
which amounts to $n-i\ge r$. So we have: $Q_{dim(A)+n-dim(E)-i}\tilde{Q}_{dim(E)+i-n}=0$ for any $\tilde{Q}$,
so $Q_{dim(A)+n-dim(E)-i}=0$ for $i>n-r$.
\epr
Another restriction comes from the symmetry of the definition \ref{DefLifting} with respect to $y_i$. $\tSi$
is invariant with respect to a subgroup $G$ of group of permutations of $y_i$ that preserves the
 flag structure. The group has the orbits:
\beq
(y_1\dots y_{dim(V_{n-1})})~~(y_{dim(V_{n-1})+1}\dots y_{dim(V_{n-2})})\dots(y_{dim(V_{2})+1}\dots y_{dim(V_{1})})
\eeq
Therefore we have:
\bcor\label{ClaimSeconConsistencCondition} {\bf (The second consistency condition)}\\
The cohomology class of~ $\tSi(x,y_1,\dots,y_r)$ is invariant under the action of the group $G$
(i.e. is a polynomial symmetric with respect to relevant subsets of $Y_i$).
\ecor
\bex
In particular for a multiple point of co-rank $r$, equation (\ref{CorankLiftings}), we have:

The cohomology class of the lifted stratum of multiple point of co-rank $r$, expressed in terms of
 $(Y_1,\dots,Y_r,X)$ and~\mbox{$Q=(d-p)X+F$}, is symmetric with respect to $(Y_1,\dots,Y_r,X)$ and the maximal
 powers of variables $(Y_1,\dots,Y_r,X)$ are not higher than $(n-r)$.
\eex
So we have rather restrictive conditions on the possible class of $\tSi$. These conditions
enable us to avoid lengthy calculations of some parameters (as will be demonstrated in section:
\ref{SecFormOfCorank}).

One could also consider the second projection: from $\tSi$ to the auxiliary space,
in definition \ref{DefLifting} this space is:
$\mP^n_x\times\prod_i\mP^n_{y_i}$.
\bex For a multiple point of co-rank $r$, the variety defined in (\ref{CorankLiftings}) is a locally trivial
fibration over
the auxiliary space $\mP^n_x\times\prod_{i=1}^r\mP^n_{y_i}$ outside the "diagonals"
(two coinciding points, three points on a line, four points in the plane etc.).
\eex
 This situation happens in general case: the projection is a locally trivial fibration, outside the cycles of jump.

\paragraph{On the possible cycles of jumps}\label{SecPossibleCyclesOfJump}
The cycles of jump were discussed in general in $\S$ \ref{SecResidualVarieties}. Here we describe the possible
cycles of jump for linear singularities.
The lifted stratum is intersection of hypersurfaces, each being defined by vanishing of a particular derivative.
As follows from the discussion in $\S$ \ref{SecDefinCondLinearQuasiHomog} the defining equations of the
hypersurfaces are of the form:
\beq\label{EqHypersurfCyclHump}
f^{(k)}|_x(\underbrace{y_{i_1},\dots,y_{i_{l_k}}}_{l_k},\hat{e}_{j_1},\dots,\hat{e}_{j_{k-l_k}})=0~~~
0\le j_1\le\dots\le j_{k-l_k}\le n
\eeq
(Here $\hat{e}_{j_1},\dots,\hat{e}_{j_{k-l_k}}$ are the vectors of the standard basis, defined in
section \ref{SecDefinCondLinearQuasiHomog}).

The cycles of jumps consist of points of $\mP^n_{x}\times\mP^n_{y_1}\times\dots\times\mP^n_{y_r}$ over which
the intersection is non-transversal. Every equation of the type (\ref{EqHypersurfCyclHump}) is linear in
the coefficients of $f$ (which are the homogeneous coordinates of the parameter space $\mPN$). Therefore the
non-transversality can happen only when the vectors $y_{i_1},\dots,y_{i_{l_k}},\hat{e}_{i_1},\dots,\hat{e}_{i_n-l_k}$ are non-generic
with respect to the vectors of other equations. More precisely:
the non-transversality can happen only when some of the vectors $y_{i_1},\dots,y_{i_{l_k}},\hat{e}_{i_1},\dots,\hat{e}_{i_n-l_k}$ of
one equation belong to the span of the vectors of other equations.

This condition can be nicely written using projections of vector space. Let $I=\{i_1,\dots,i_k\}$ be an arbitrary
(non-empty) subset of $\{0,1,\dots,n\}$. Represent the point of a projective space by its homogeneous coordinates:
$y=(z_0,\dots,z_n)\in\mP^n_y$. The projection is defined by:
\beq
\pi_I(y):=(z_{i_1},\dots,z_{i_k})
\eeq
Note that $\pi_I(y)$ is defined up to a scalar multiplication and can have all the entries zero. Immediate
check shows that the above condition of non-transversality corresponds to one of the following:
\bei
\item $\pi_I(y)\in$Span$(\pi_I(x),\pi_I(y_{i_1}),\dots,\pi_I(y_{i_k}))$
\item $(\pi_I(x),\pi_I(y_{i_1}),\dots,\pi_I(y_{i_k}))$ are linearly dependent.
\eei
(In particular if all the entries of $\pi_I(y)$ are zero, both conditions are trivially satisfied.)
Note that the second condition is the closure of the first. Summarizing:
\bprop\label{ClaimPossibleCyclesOfJump}
The possible cycles of jump in the auxiliary space $Aux=\mP^n_x\times\mP^n_{y_1}\times\dots\times\mP^n_{y_k}$
are of the form: $\Bigg((\pi_I(x),\{\pi_I(y_j\}_{j\in J})$ are linearly dependent for some
$I\subseteq\{0,\dots,n\},~~J\subseteq\{1,\dots,k\}\Bigg)$
\eprop
\beR
While the set $I$ is nonempty, $J$ can be empty. In this case the cycle of jump is defined as  $\pi_I(x)=(0\dots0)$.
\eeR
\paragraph{On the adjacency}\label{SecOnTheAdjacency}
To each equisingular stratum some strata of higher singularities are adjacent (i.e.the strata of higher
singularities are included in the closure of the given stratum).
For example $\Sigma_{D_{k+1}}\subset\overline\Sigma_{A_k}$.

We constantly use the codimension one adjacency, i.e. when $\Si_{\mD'}$ is a divisor in $\Si_\mD$
(or the same for the lifted versions, cf. $\S$ \ref{SecLiftings}).
In other words, these are the strata that can be reached by just one degeneration.
Many tables of adjacencies are given in \cite{AVGL}. In each particular case the adjacency can be
checked by the analysis of Newton diagram of the singularity or of the defining ideal of the singular germ.

The adjacency can depend on moduli \cite{Pham}. However, if one chooses moduli generically (and we always
do that), the adjacency is completely fixed by the topological type.

A more important feature is: the set of relevant adjacent strata can depend on the dimensionality of
the ambient space. As an example consider the enumeration of $A_4$ (cf. \ref{SecFurtherCalculations}).
It is degenerated by increasing the corank of
the quadratic form. Then:\\
$\bullet n=2$, the case of curves $A_4\rightarrow D_5$\\
$\bullet n\geq3$, the case of hypersurfaces, $A_4\rightarrow D_5\cup P_8$
\\
As will be shown in $\S$ \ref{SecHigherSing}, in the second case both $D_5$ and $P_8$ are relevant.

In our approach we enumerate separately for each fixed dimension $n$.
By universality, the final answer is given by a unique polynomial (in Cherna
classes that depend on $n$) of a known degree. Therefore it suffices to calculate just for a
few needed dimensions and by universality we get the complete answer.

\subsubsection{The ideology of degenerating process}\label{SecIdeologyOfDegenerations}
Here we prove the main theorems stated in Introduction (\ref{TheoremEnumLinearSings} and \ref{TheoremEnumNonLinearSings}). We must prove
 the three statements:\\
$\bullet$the degenerating step is always possible and the degeneration is invertible.\\
$\bullet$for linear singularities we achieve the cohomology class of the lifted stratum in a finite number of steps
and all the intermediate types are linear.\\
$\bullet$for non-linear singularities we reduce the problem to enumeration of linear ones (in a finite number of steps).
\paragraph{The degenerating step}consists of intersection of a lifted stratum with a hypersurface and subtraction
of the residual pieces. In more details, let $\tSi$ be the initial lifted stratum and $\tSi_{degen}$
the stratum we want to reach.  The degenerating hypersurface is defined by the equation
$f^{(p)}(y_{i_1},\dots,y_{i_k})_{j_1\dots j_{p-k}}=0$. For linear singularities this is just the derivative
corresponding to a point under the Newton diagram. For non-linear singularities (as is explained
in $\S$ \ref{SecDefinCondNonLinearSing}) the derivative
corresponds to a point on the Newton diagram with the minimal distance to the origin. Such a point
can be non-unique, in this case we use monomial order from $\S$ \ref{SecOnTheMonomialOrder}.

The intersection is in
general non-transversal and the resulting variety is reducible (containing residual pieces in addition
to the needed degenerated stratum).
\beq
\tSi\cap \{f^{(p)}(y_{i_1},\dots,y_{i_k})_{j_1\dots j_{p-k}}=0\}=\tSi_{degen}\cup\{\mbox{Residual pieces}\}
\eeq
The non-transversality happens over cycles of jump (described in $\S$ \ref{SecPossibleCyclesOfJump}).
The procedure to calculate the cohomology classes of the residual pieces is explained in section
\ref{SecCohomologyClassOfRestrictionFibration}. Note that in the above equation the degenerated stratum
 $\tSi_{degen}$ can be reducible and non-reduced (but is always pure dimensional), this happens for non-linear singularities. The multiplicity
 of the intersection of $\tSi$ with the degenerating hypersurface along $\tSi_{degen}$ is calculated
 in the classical way (section \ref{SecMultiplicityOfIntersection}).

Summarizing, the degenerating step produces the equation in cohomology:
\beq
[\tSi][f^{(p)}(y_{i_1},\dots,y_{i_k})_{j_1\dots j_{p-k}}=0]=[\tSi_{degen}]+[\mbox{Residual pieces}]
\eeq
By $\S$ \ref{SecInvertibilityofDegeneration} the degeneration is invertible, so the equation enables the calculation of
either $[\tSi]$ or $[\tSi_{degen}]$.
\paragraph{Degenerations for linear singularities (from simple to complicated).}
The defining set of conditions for linear singularities was described in $\S$ \ref{SecDefiningConditionsGeneral}.
Starting from the stratum of ordinary multiple point (of the relevant multiplicity), lifted to the space $Aux\times\mPN$,
we apply the defining conditions one-by-one. Each condition means
the absence of a particular monomial, arrange the conditions by the monomial order
(defined in \ref{SecOnTheMonomialOrder}). This guarantees that at each step we get a Newton-non-degenerate type.

The degenerating step was described above. After a restricted number
of degenerations (not bigger than the number of points under the Newton diagram)
we arrive at the lifted stratum of the needed singularity.

At each intermediate step of the process the singularity type is linear. Indeed by the criterion
 \ref{ClaimCriterionLinearity} all the initial and final slopes are bounded in the interval $[\frac{1}{2},2]$.
 And in the process of degeneration the slopes change monotonically.
\beR
Usually it is simpler to start not from the stratum of ordinary multiple point, but from a stratum of a higher (linear) singularity
to which the given singularity is adjacent and for which the enumeration problem is already solve. For example
to enumerate the tacnode ($A_3$) we start from the cusp ($A_2$).
\eeR
\paragraph{Degenerations for non-linear singularities (from complicated to simple)}
As was explained in $\S$ \ref{SecSimplifyingDegenerations} the original non-linear type is
degenerated to a combination of linear ones.

The goal of degenerating process is to convert the defining non-linear equations to monomial ones (i.e. of the
form $\vec{z}^{\vec{m}}=0$). For this, at each step of the process, we consider a (non-vanishing)
monomial of the lowest monomial order (section \ref{SecOnTheMonomialOrder}).

Suppose the multiplicity of the singularity is $p$ while the order of determinacy is $k$.
The process goes by first demanding that the derivatives $f^{(p)}_{i_1\dots i_p}$ vanish (the conditions are
applied one-by-one, at each step we have just the degeneration by a hypersurface). Then one arrives at
the singularity of multiplicity $(p+1)$ and order of determinacy $k$. If the so obtained singularity
(or a collection of singularities) is still non-linear, the process is continued.
In the simplest case (Newton-non-degenerate singularity), once we have $k\le 2p$, we necessarily have
a collection of linear singularities (by corollary \ref{CorollarLinearSingMultPOrderDeterm2p}).
In the worst case (Newton-degenerate singularity) one continues up to the ordinary point of multiplicity $k$.

This process transforms a non-linear singularity to a collection of linear singularities whose
enumeration was described above. As all the degenerations are invertible, this solves the enumerative
problem for non-linear singularities.
\subsection{The simplest examples}\label{SecExamples}
We consider here some simplest typical examples to illustrate the method.
First we consider the case of cusp. Having enumerated the cusp, one
enumerates the tacnode by just one additional degeneration.
\subsubsection{Quadratic forms of co-rank 1 (cuspidal hypersurfaces)}
Here we consider singularity with the normal form $\sum_{i=2}^nz^2_i+z^3_1$.
The lifted stratum was defined in example \ref{ExampleNumerousLiftings}:
\beq
\tSi_{A_2}(x,y)=\overline{\{(x,y,f)|~(x\ne y),~f|_x^{(2)}(x)=0=f|_x^{(2)}(y)\}}\subset\mP^n_x\times\mP^n_y\times\mPN
\eeq
(here the tensor of second derivatives is calculated at the point $x$, therefore: $f|_x^{(2)}(x)\sim f|_x^{(1)}$).

We want to represent $\tSi_{A_2}(x,y)$ as a (possibly) transversal intersection of hypersurfaces.
The $(n+1)$ conditions: $f^{(2)}(x)=0$ are transversal. Suppose we add to them one
condition: $(f^{(2)}(y))_0=0$. Then, the non-transversality occurs over two cycles of jump in $Aux=\mP^n_x\times\mP^n_y$:
\bei
\item{$x=y$} (co-dimension $n$). The jump in the dimension of the fiber is 1.
\item{$x=(1,0,\dots,0)$} (co-dimension $n$). In this case, since $f^{(2)}(x)=0$, we already know that
the form $f^{(2)}(y)$ annihilates $x$ (i.e. $f^{(2)}(y)_0=0$). The jump in the dimension of the fiber is 1.
\eei
In both cases the dimension of the jump of fiber is less than the codimension of the cycle of jump,
so the resulting variety is irreducible. So, for the cohomology classes we have:
\beq
[\overline{f^{(2)}(x)=0,f^{(2)}(y)_0=0,~x\neq y}]=[f^{(2)}(x)=0]\times[f^{(2)}(y)_0=0]
\eeq
We continue to degenerate by the conditions in such a way up to $(f^{(2)}(y)_{n-1}=0)$. At this point
for generic (non-coinciding) $x,y$ we have:
$f^{(2)}(x)=0=f^{(2)}(y)$.
\\
\parbox{13cm}{Additional pieces that arose are:
\\$\bullet$ over $x=y$, the condition of codimension $n$. Over this diagonal the jump in the dimension of fibers is  $n$.
\\$\bullet$ over $x=(*\dots*,0)$, the condition of codimension 1. Over this subvariety the jump in the
dimension of fibers is $1$.\\
On the picture the result of the intersection is shown in a "log-scale", i.e. the dimensionality is
transformed to the relative hight.}
\begin{picture}(0,0)(0,20)
\ellipse{0}{15}{30}{0}{120}{20}{80}{40}
\put(10,2){\line(3,1){100}}\put(10,42){\line(3,1){100}}\put(15,15){$x=y$}
\put(10,2){\line(0,1){40}} \put(107,34){\vector(0,1){40}}\put(107,74){\vector(0,-1){40}}  \put(115,50){$n$}

\put(110,10){\line(-3,1){40}}\put(110,20){\line(-3,1){40}}\multiput(65,24)(-6,2){6}{.}\multiput(65,34)(-6,2){6}{.}
\put(107,13){$\updownarrow$} \put(110,0){$x_n=0$}\put(112,12){$1$}
\put(45,-15){$\mP^n_x\times\mP^n_y=Aux$}
\end{picture}
\\
So, the resulting variety is reducible, containing residual pieces over the cycles of jump:
\beq
\Big\{f|_x^{(2)}(x)=0=f|_x^{(2)}(y)\Big\}=\tSi_{A_2}(x,y)\cup \Big\{Res_{x=y}\Big\}\cup \Big\{Res_{z_n=0}\Big\}
\eeq
(By direct check we obtain that the multiplicity in both cases is 1.)
To remove the contributions from additional pieces we should calculate the classes of the restrictions to the cycles of jump.
The method of calculation of the cohomology classes of residual pieces is described in section
\ref{SecCohomologyClassOfRestrictionFibration}.
\\
$\bullet$ The residual piece over the diagonal can be described explicitly in a simple manner. As all the
intersections over the points of the diagonal are nontransversal, a point over the diagonal corresponds to a nodal
hypersurface. So, the residual piece is $Res_{x=y}=\{x=y\}\cap\Big\{\tSi_{A_1}(x)\times\mP^n_y\Big\}$,
and its cohomology class $[Res_{x=y}]$=$[x=y]\times$$[\tSi_{A_1}(x)]$.
\\
$\bullet$  Over the generic point of the cycle $z_n=0$ (i.e. the point for which $z_0\dots. z_{n-1}\ne 0$) all
the intersections, except for the last were transversal. Therefore we have:
\beq
\Big\{f|_x^{(2)}(x)=0,~~z_n=0\Big\}\cap^{n-2}_{i=0}\Big\{f|_x^{(2)}(y)_i=0\Big\}=Res_{z_n=0}\cup Res_{z_n=0=z_{n-1}}
\eeq
where $Res_{z_n=0=z_{n-1}}$ is a "secondary" residual piece. Its description is the same as that of $Res_{z_n=0}$
and so one has a recursion:
\beq
\Big\{f|_x^{(2)}(x)=0,~~z_n=0\Big\}\cap^{j-2}_{i=0}\Big\{f|_x^{(2)}(y)_i=0\Big\}=Res_{z_n=0=\dots=z_j}
\cup Res_{z_n=0=\dots=z_{j-1}}
\eeq
After completion of the recursion we get the formula:
\beq
[\tSi_{A_2}(x,y)]=[f^{(2)}(x)=0]\Bigg(\sum_{i=0}^{n-1}(-1)^i[\prod_{j=n+1-i}^n\{x_j=0\}]
[\prod_{j=0}^{n-1-i}\{f^{(2)}(y)_j=0\}]-[x=y]\Bigg)
\eeq
Or, in terms of cohomology classes:
\beq\label{CuspCoefficient}
[\tSi_{A_2}(x,y)]=(Q+X)^{n+1}\left(\sum^n_{i=0}(-1)^i(Q+Y)^{n-i}X^i-\sum^n_{i=0}X^iY^{n-i}\right),~~Q=(d-2)X+F
\eeq
Here $(X,Y,F)$ are the generators of the corresponding cohomology rings. Note, that this polynomial
(if written in variables $(Q,X,Y)$)is
symmetric in $X,Y$ (eventhough it was obtained in a very non-symmetric way) as it should be. Also, the terms $X^n,Y^n$
cancel in the polynomial.

To obtain the degree of the variety we should extract the coefficient of $X^nY^{n-1}$ (after the substitution: $Q=(d-2)X+F$).
We get:
\beq\label{cuspdeg}
\mbox{deg}(\Sigma_{A_2})=(d-1)^{n-1}(d-2)\frac{n(n+1)(n+2)}{2}
\eeq
(which of course coincides with the result obtained by P.Aluffi in \cite{Alufi}).
\subsubsection{The use of adjacency: tacnodal hypersurfaces}\label{SecTacnodalHypersurf}
The tacnodal singularity ($A_3$) has the normal form $f=z^4_1+\sum_{i=2}^nx_i^2$.
The corresponding lifted variety was defined in example \ref{ExampleNumerousLiftings}.
We represent the tacnode as a degeneration of the cusp. Correspondingly, we think of the (lifted)
stratum of tacnodal hypersurfaces as a subvariety of the cuspidal stratum:
\beq
\tSi_{A_3}(x,y)=\overline{\{(x,y,f)\in\tSi_{A_2}(x,y),~x\ne y~~f|_x^{(3)}(y,y,y)=0\}}\subset\tSi_{A_2}(x,y)\subset\mP^n_x\times\mP^n_y\times\mPN
\eeq
For generic $x,y$ the intersection is $\tSi_{A_2}(x,y)\cap\{f^{(3)}(y,y,y)=0\}$ transversal.
The possible non-transversality can occur only when $x=y$.
So for the cohomology classes we have:
\beq
[\tSi_{A_3}(x,y)]=[\tSi_{A_2}(x,y)][f^{(3)}(y,y,y)=0]-[\mbox{Residual piece over }x=y]
\eeq
The method to calculate the cohomology class of the residual piece over $(x=y)$, is given in Appendix B. By
the Corollary \ref{CorCohomologyClassResidualPieceOverDiagonal} we have:
\beq
[\mbox{Residual piece over }x=y]=[x=y]\frac{1}{(n-1)!}\frac{\di^{n-1}[\tSi_{A_2}(x,y)]}{\di^{n-1}Y}
\eeq
In this simple case the residual piece can be also easily described explicitly: every point of
it corresponds to a hypersurface with
a cusp at a given point, but with arbitrary tangent line. So the variety over the diagonal ($x=y$), is just
$\tSi_{A_2}(x)$ taken with some multiplicity. The multiplicity can be computed in two ways:
\li{directly,} as the degree of tangency of the nontransversal intersection.
\li{via consistency condition} (as was explained in $\S$ \ref{SecLiftings}). Since the lifted stratum
$\tSi_{A_3}$ is a fibration over $\Sigma_{A_3}$ with fiber
$\mP^1_y$, the corresponding cohomology class should not include terms with $Y^n$).

Both methods give the multiplicity 3. Thus the cohomology class of the cuspidal stratum is:
\beq\ber
[\tSi_{A_3}(x,y)]=[\tSi_{A_2}(x,y)][f^{(3)}(y,y,y)=0]-3[x=y][\tSi_{A_2}(x)]=
\\
=(Q+X)^{n+1}\Bigg(\Big(\sum^n_{i=0}(Q+Y)^{n-i}(-X)^i-\sum^n_{i=0}X^iY^{n-i}\Big)(Q+3Y-X)-
3(nQ-2X)\sum^n_{i=0}X^iY^{n-i}\Bigg)
\eer\eeq
Here, as always $Q=(d-2)X+F$. Note again, that (in full accordance with the fibration conditions)
the answer (if written in terms of $X,Y,Q$) is symmetric in $X,Y$ and no terms with $X^n$ or $Y^n$ appear.
Finally, the degree is:
\beq
[\Sigma_{A_3}]={n+2\choose{3}}(d-1)^{n-2}\Bigg(\frac{(3n-1)(n+3)}{2}(d-2)^2+2(n-1)(d-2)-4\Bigg)
\eeq
which for $n=2$, curves, coincides with the result of Aluffi \cite{Alufi}.
\subsection{Further calculations}\label{SecFurtherCalculations}
We start from homogeneous forms of some (co-)rank (section \ref{SecFormOfCorank}).
For quadratic forms ($p=2$) the rank fixes the degeneracy class completely. For higher forms one can impose various
additional degeneracy conditions (e.g. for $p=3$ the form: $z_1^2z_2+\sum_{i=3}^nz_i^3$ is of full rank, but
the corresponding singularity is not an ordinary triple point). We consider examples of such singularities in
section \ref{SecHigherSing}.

Having calculated the classes of the lifted strata for $2-$forms of some rank, we can start enumeration of
other singularities. This is done by further degenerations. For example, the tacnode $A_3$ is the cusp with some
degeneracy of the tensor of order-3 derivatives (it was enumerated in $\S$ \ref{SecExamples}).
We consider here the simplest examples: $A_4,D_5,E_6$.
\subsubsection{The forms of co-rank$\geq1$}\label{SecFormOfCorank}
We first recall:
\bed
The homogeneous symmetric form of order $p$, in $n$ variables, $\Omega^{p}(z_1,\dots,z_n)$, is called {\bf of rank $(n-r)$}
(of co-rank $r$)
if by linear transformation of $GL(n)$ in the space of variables ($\mC^n=\{z_1,\dots,z_n\}$) it can be brought to
a homogeneous form in $n-r$ variables: $\tilde\Omega^{(p)}(z_1,\dots,z_{n-r})$.
\eed
The collections of such forms are natural generalizations of classical determinantal varieties (symmetric matrices of a given co-rank)
\cite[chapter 14.3]{Ful}. To emphasize the co-rank of the form we often assign it as a subscript: $\Omega^{p}_r$.

There are (at least) two approaches to calculate the cohomology classes of the stratum  of a given co-rank forms:
\li{Start from the ordinary multiple point (it corresponds to the non-degenerate form) treated before}. Apply the
degenerating conditions (one-by-one) to get the form of the needed co-rank.
At each step it is necessary to remove the residual pieces. This approach works well for the forms
of low co-rank.
\li{Degenerate the given form to a form of rank 2 or 1}. (The forms of rank$=2,1$ are particular cases
of reducible forms, their enumeration is immediate and is treated in $\S$ \ref{SecReducibleForms})
Then we will have equation in cohomology of the form:
\beq
[\Omega^{(p)}_r][\mbox{degenerating cycle}]=[\Omega^{(p)}_2]+[\mbox{Residual variety}]
\eeq
And from this equation the class $[\Omega^{(p)}_r]$ is restored uniquely.

We describe here the first approach. As the computations are extremely involved
we solve explicitly only the case of the quadratic forms of arbitrary corank.

The lifted stratum was defined in example \ref{ExampleNumerousLiftings}:
\beq\label{CuspVector}
\tSi^n_r(x,(y_i)_{i=1}^r)=\overline{\Big\{(x,\{y_i\}_{i=1}^r,f)\Big|
\Big(\!\!\ber
x,y_1,\dots,y_r~\rm{are}\\\rm{linearly~independent}\eer\!\!\Big),
\Big(\!\!\ber~~~f^{(2)}(x)=0\\(f^{(2)}(y_i)=0)_{i=1}^r\eer\!\!\Big)\Big\}}\subset\mP^n_x\times\prod_{i=1}^r\mP^n_{y_i}\times\mPN
\eeq
As was explained in $\S$ \ref{SecLiftings}, the cohomology class of $\tSi^n_r$ is a polynomial in
$(X,Y_1,\dots,Y_r$, $Q=(d-2)X+F)$, symmetric in $(X,Y_1,\dots,Y_r)$ and does not contain powers of $Y_i$,
greater than $(n-r)$. The cohomology class of $\Sigma^n_r$ is just the coefficient of the monomial
$Y_1^{n-r}\dots Y_k^{n-r}X^n$ in the cohomology class of $\tSi^n_r$.

The enumeration of singularities with quadratic form of co-rank $r$ is completed by the claim:
\bel
The cohomology class of the lifted variety $\tSi^n_r(x,(y_i)_{i=1}^r)$ can be calculated by
successive degenerations starting from $\tSi^n_0(x)$
(the nodal hypersurfaces). In particular, the cohomology class of the minimal lifting is:
\beq\ber
[\tSi^n_r(x)]=C_{n,r}Q^{r\choose{2}}
\sum^r_{i=0}\frac{{n-i\choose{r-i}}{r\choose{i}} }{{2r\choose{r+i}}}Q^{r-i}(-X)^i~~~~~~Q=F+(d-2)X\\
C_{n,1}=2,~~C_{n,2}=2{n+1\choose{1}},~~C_{n,3}=2{n+2\choose{3}},~~C_{n,4}=2{n+3\choose{5}}\frac{n+1}{3}\\
C_{n,n}=2{2n\choose{n}},~~C_{n,n-1}=\frac{2^r{2r\choose{r}}}{n},~~
C_{n,n-2}=\frac{{2(r+1)\choose{r+1}}{2r\choose{r}}}{{r+2\choose{2}}}
\eer\eeq
\eel
Note that here we give the constant for some specific values of $(n,r)$ only. This is due to computer limitations:
each time we calculate for a specific value of $r$ and $n$. So, for every specific $n,r$ one can get the
answer (provided the computer is strong enough), but it is not clear how to combine these values
into one nice expression.
\bpr The case $r=1$ corresponds to cuspidal hypersurfaces and was considered in $\S$ \ref{SecExamples}.
The general case is done recursively. Suppose we have obtained the cohomology class of
$\tSi^n_r(x,y_1,\dots,y_r)$,
as in equation (\ref{CuspCoefficient}). Intersect the variety $\tSi^n_r(x,y_1,\dots,y_r)$
with $(n-r)$ hypersurfaces: $f^{(2)}(y_{r+1})_0=0,\dots,f^{(2)}(y_{r+1})_{n-r-1}=0$. The possible (significant)
non-transversality can occur
in two cases, either: $y_{r+1}\in$Span($x,y_1,\dots,y_r$), or
$(*,\dots,*,\underbrace{0,\dots,0}_k)$$\in$Span($x,y_1,\dots,y_r$). In both cases one continues as in the case of $r=1$.
In such a way we get the cohomology class:
\beq
[\tSi^n_{r+1}]=[\tSi^n_r]\Big(\sum_{i=0}^{n-r}(Q+Y_{r+1})^{n-r-i}(-1)^i
\hspace{-0.3cm}\sum_{j_1+\dots+j_r\leq
i}\hspace{-0.3cm}Y_1^{j_1}\dots Y_r^{j_r}X^{i-(j_1+\dots+j_r)}-
\hspace{-0.9cm}
\sum_{i_1+\dots+i_{r+1}=n-r-j}\hspace{-0.9cm}X^{j}Y_1^{i_1}\dots
Y_{r+1}^{i_{r+1}}\Big)+\rm{residual~terms}
\eeq
Here the residual terms correspond to varieties which occur over the diagonals: $x=y_1$,  $y_2\in$span($x,y_1$),
 \dots,
$y_r\in$span($x,y_1,\dots,y_{r-1}$). These residual pieces can be calculated by classical intersection theory
(as explained in $\S$ \ref{SecCohomologyClassOfRestrictionFibration}). However (as happens in all other cases)
the classes
are actually completely fixed by the consistency conditions (lemma \ref{ClaimFirstConsistencyCondition} and
\ref{ClaimSeconConsistencCondition}). In this case they read:

{\it The final expression (polynomial in $Q,X,Y_1,\dots,Y_{k+1}$) should be symmetric in $(X,Y_1,\dots,Y_{k+1})$ and
should not contain the powers of $X,Y_1,\dots,Y_{k+1}$ that are greater than $n-k-1$}.

The explicit calculations are extremely complicated and can be done by computer only. The cohomology classes of
the lifted strata are awkward polynomials in many variables.

For example the cohomology class in the co-rank 2 case ($D_4$ singularity) is:
\beq\ber
[\tSi^2_{2}(x,y_1,y_2)]=\tinyM(Q+X)^{n+1}
\sum_{i=0}^n(Q+Y_1)^{n-i}(-X)^i\Big(\sum_{j=0}^{n-1}(Q+Y_2)^{n-1-j}(-1)^j
\sum_kX^{j-k}Y_1^k-\sum_{j,k}X^jY_1^kY_2^{n-1-j-k}\Big)
-\\-\tinyM
(Q+X)^{n+1}\sum_{i=0}^nX^iY_1^{n-i}\sum^n_{i=1}\Big(\sum^{n-i}_{j=0}{n-j\choose{i}}Q^{n-i-j}(-X)^j-X^{n-i}\Big)
\Big(
\sum_{j=0}^{i-1}(-1)^{n-j-1}(Q+Y_2)^{j}X^{i-j-1}-\sum_{j=0}^{i-1}X^{i-1-j}Y_2^{j}\Big)
\eer\eeq
Extracting the coefficients of $Y_1^{n-r-1}\dots Y_{k+1}^{n-r-1}$ we get the classes of the strata $\tSi^n_r(x)$.
\epr
\subsubsection{Some linear singularities}\label{SecHigherSing}
\paragraph{Hypersurfaces with a $D_5$ point,} the normal form: $f=z^4_1+z^2_2z_1+\sum_{i=3}^nz_i^2$.
The lifted stratum was defined in example \ref{ExampleNumerousLiftings}.
We represent the stratum $\tSi_{D_5}$ as a subvariety of $\tSi_{D_4}$:
\beq
\tSi_{D_5}(x,y_1,y_2)=\overline{\{(x,y_1,y_2,f)\in\tSi_{D_4}(x,y_1,y_2),~~f|_x^{(3)}(y_1,y_1,y_1)=0=f|_x^{(3)}(y_1,y_1,y_2)\}}
\subset\tSi_{D_4}(x,y_1,y_2)
\eeq
For generic $x,y_1,y_2$ the intersection is transversal. The non-transversality occurs over diagonals:
\mbox{$y_2\in$Span$(x,y_1)$}
or $x=y_1$. Note, that the first variety is non-closed. We approximate it by the variety: ($x,y_1,y_2$ are
linearly dependent),
the two varieties coincide for $x\neq y_1$. Correspondingly, over $x=y_1$ we have additional (secondary) residual piece.
Thus the cohomology class is:
\beq
[\tSi_{D_5}(x,y_1,y_2)]=\Big(\ber[\tSi_{D_4}(x,y_1,y_2)][f^{(3)}(y_1,y_1,y_1)=0][f^{(3)}(y_1,y_1,y_2)=0]-\\-
[x=y_1]A(X,Y_2,Q)-[\rm{rk}\Big(\begin{smallmatrix}x\\y_1\\y_2\end{smallmatrix}\Big)<3]B(X,Y_1,Y_2,Q)
\eer\Big)
\eeq
Here $A,B$ are some (homogeneous) polynomials in the generators of the cohomology ring. By the identity in the
cohomology ring:  $(X-Y)\sum_{i=0}^n X^iY^{n-i}=X^{n+1}-Y^{n+1}\equiv0$, we can assume that $A$ does not depend on $Y$.
The only additional condition on $A,B$ is the consistency condition from $\S$ \ref{SecLiftings}:

{\it The cohomology class $[\tSi_{D_5}(x,y_1,y_2)]$ should not contain monomials with $Y^n_1,Y^{n-1}_2,Y^n_2$.}

As always, this condition itself fixes the polynomials completely. The final cohomology class is given in
Appendix.
\paragraph{Hypersurfaces with an $E_6$ point,}the normal form: $f=z^4_1+z^3_2+\sum_{i=3}^nz_i^2$.
We represent $E_6$ as a degeneration of $D_5$:
\beq
\tSi_{E_6}(x,y_1,y_2)=\overline{\{(x,y_1,y_2,f)\in\tSi_{D_5}(x,y_1,y_2),~~f|_x^{(3)}(y_1,y_2,y_2)=0\}}
\subset\tSi_{D_5}(x,y_1,y_2)
\eeq
Again, instead of describing the residual varieties explicitly, we use the consistency conditions, which completely
fix the class. The final answer is in Appendix.
\subsubsection{Some non-linear singularities}
\paragraph{The $A_4$ case}
Here we consider the simplest non-linear case. By linear transformation, the singularity germ can be brought
to the Newton diagram of $A_3$: $f=\sum_{i=2}^n\alpha_iz_i^2+z_1^2\sum_{i=2}^n\beta_iz_i+\gamma z^4_1+\dots$

To achieve the Newton diagram of $A_4$ we must do the non-linear shift: $z_i\rightarrow z_i-\frac{\beta_i}{2\al_i}z^2_1$
(to get rid of the monomials $z^4_1,z^2_1z_i$). Elimination gives $\gamma=\sum\frac{\beta^2_i}{4\alpha_i}$.

Therefore, degenerating $\alpha_2=0$ we get: $\beta_2=0$ or $\prod_{i\ne 2}\alpha_i=0$, corresponding to adjacency
$\bar\Si_{A_4}\supset\Si_{D_5},\Si_{P_8}$

In this way we obtain the cohomology class of $\tSi_{A_4}(x,y)$. So, we get the equation for
cohomology classes:
\beq
[\tSi_{A_4}(x,y_1)][degeneration]=2[\tSi_{D_5}(x,y_1)]+2[\tSi_{P_8}(x,y_1,y_2)]+[residual~~piece]
\eeq
The final result is in Appendix.
\paragraph{The $D_6$ case}
By linear transformation, the singularity germ can be brought to the Newton diagram of $D_5$:
$f=\sum_{i=3}^n\alpha_iz_i^2+z_1^2\sum_{i=3}^n\beta_iz_i+\gamma z^4_1+z_1z^2_2\dots$
To achieve the Newton diagram of $D_6$ we must do the non-linear shift: $z_i\rightarrow z_i+\delta_iz^2_1$
(to get rid of the monomials $z^4_1,z^2_1z_i$). Elimination gives: $\gamma=\sum\frac{\beta^2_i}{4\alpha_i}$.
We degenerate in the same way as in the $A_4$ case and get:
\beq
[\tSi_{D_6}(x,y_1,y_2)][degeneration]=2[\tSi_{P_8}(x,y_1,y_2,y_3)_{degenerated}]+
[corank~4]+[residual~~piece]
\eeq
We omit the calculations.
\appendix
\section{Singularities with reducible jets}\label{SecReducibleForms}
Here we consider singular polynomials whose low order jet is reducible. The simplest such case is that of reducible form:
\beq\label{ReducibForm}
\Omega^{(p)}=\prod_{i=1}^k\Bigg(\Omega^{(p_i)}_i\Bigg)^{r_i},~~~\sum_{i=1}^k r_ip_i=p
\eeq
(here the homogeneous forms $\Omega^{(p_i)}_i$ are irreducible, though they can be degenerate and mutually non-generic).
These singularity types are of high codimension, therefore extremely rare, nevertheless
they deserve some attention, being sometimes the final goal of degenerating process.

Since reducibility is in general not invariant under the topological transformations we (in general)
cannot define the corresponding stratum as a topological one. We define the stratum as the collection
of hypersurfaces that can be brought (by locally analytic transformation) to  {\it of a given form}.

This stratum is included into the topological stratum. On the other side it usually contains families of
 analytical strata, since by the Newton diagram we do not specify the moduli.

In case of curves ($n=2$) every singularity of multiplicity $p$ has reducible $p-$jet
(being homogeneous polynomial of two variables). So the corresponding stratum coincides with
the topological equisingular stratum.
Some singularities with reducible $p-$form are: $A_1,A_2,D_4$,$E_6,X_9$,$Z_{11},W_{12}\dots$.

For surfaces ($n=3$) some singularities of this type are: $A_2,D_4,P_8,S_{11}$,
$U_{12},T_{3,4,4},T_{4,4,4}$,$V_{1,0},V'_1\dots$.

Another example of reducible p-form (for any $n$) is the p-form of rank 2 or 1.

For singularities with reducible jets, the lifted stratum can be explicitly defined by conditions of
a very standard type: proportionality of tensors.
\subsubsection{Reducible homogeneous forms}\label{SecDegenerateHomogenForms}
 We consider here the case of mutually generic forms $\Omega^{(p_i)}_i$ in (\ref{ReducibForm}).
Every such form defines a hypersurface, by the mutual generality of the forms the hypersurfaces
intersect in a generic way, however each hypersurface can be singular. We restrict to the case of ordinary
multiple point of maximal multiplicity (so for each hypersurface the condition is: $\Omega^{(p_i)}_i(x)=0$).
The stratum of hypersurfaces with this type of reducibility is defined as:
\beq\label{VarietyReducForm}
\tSi\Big(x,(\Omega^{(p_i)}_i)_{i=1}^k\Big)=\Bigg\{\Big(x,(\Omega^{(p_i)}_i)_{i=1}^k,f\Big)\Big|
~f^{(p)}\sim\rm{SYM}\Big((\Omega^{(p_1)}_1)^{r_1},\dots,(\Omega^{(p_k)}_k)^{r_k}\Big),
~(\Omega^{(p_i)}_i(x)=0)_{i=1}^k\Bigg\}
\eeq
Here SYM means symmetrization of indices.
Note that the sets of defining conditions are mutually transversal, e.g. $f$ appears in the first proportionality
condition only. Therefore the cohomology class of the lifted stratum is just the product of classes of conditions:
\beq
[\tSi\Big(x,(\Omega^{(p_i)}_i)_{i=1}^k\Big)]=[f^{(p)}\sim\rm{SYM}
\Big((\Omega^{(p_1)}_1)^{r_1},\dots,(\Omega^{(p_k)}_k)^{r_k}\Big)]\prod_i [\Omega^{(p_i)}_i(x)=0]
\eeq
The condition of proportionality of two tensors is considered in $\S$ \ref{SecCohomologyClassesOfCyclesofJump}
(equation (\ref{EqProportionalityTwoTensors})). In terms
of the cohomology ring generators of the ambient space ($X,Q=(d-p)X+F,\Omega_i$) we have:
\beq\label{EqClassesForReducibleForms}
[\tSi\Big(x,(\Omega^{(p_i)}_i)_{i=1}^k\Big)]=(\Omega_1+X)^{p_1-1+n\choose{n}}\dots(\Omega_k+X)^{p_k-1+n\choose{n}}
\sum_{i=0}^{{p+n\choose{n}}-1}Q^i(r_1\Omega_1+\dots+r_k\Omega_k)^{{p+n\choose{n}}-1-i}
\eeq
Note, that depending on the singularity type, the projection: $\tSi\rightarrow\Sigma$ can be not 1:1. The
permutation group of forms of equal multiplicity and degeneracy acts on fibers. Thus, to obtain the cohomology class of
$\Sigma$ one should divide the corresponding coefficient by the order of this group: $|Aut|$.
\subsubsection{Singularities with reducible jets}\label{SecReducibleJets}
Here we consider singularities of the type:
\beq\label{ReducibJet}
f=\prod_{i=1}^k f_i(z_1,\dots,z_n)+higher~order~terms~~~~~\rm{deg}(f_i)=p_i~~\sum_{i=1}^kp_i=p
\eeq
The polynomials $f_i$ are (non-homogeneous) of fixed degrees. We assume that the singularities $f_i=0$ are linear,
in particular they satisfy: $2\times multiplicity\ge p_i$.

In particular, we can assume that all the hypersurfaces $\{f_i=0\}$ pass through the origin. Introducing factors that do not vanish at the
origin leads to  hypersurfaces with flexes, (the property which is not invariant under local diffeomorphism/homeomorphism).

 If the hypersurface ($f_i=0$) is smooth and generic with respect to other hypersurfaces
 (e.g. all the normals are in general position, intersection is along generic subvarieties etc.) then deg($f_i$)$=p_i=1$.

The procedure of enumeration is as in $\S$ \ref{SecReducibleForms}: the problem is reduced to enumeration
of particular singularities ($f_i$), if some of the hypersurfaces are in a mutually special position this should be also taken
into account.

We consider some typical situations:
\li{\bf Mutually generic smooth hypersurfaces.} As was explained above, in this case  all the degrees
are necessarily equal to 1. So, all the factors are linear, this case was treated in section
\ref{SecReducibleForms}.
\li{\bf Mutually generic singular hypersurfaces}. In this case every singular hypersurface is treated separately,
then the results are combined. The simplest case is:
\beq
f=g\prod_{i=1}^k f_i+higher~order~terms
\eeq
Here $f_i$ are homogeneous polynomials, while $g$ is not necessarily homogeneous, with the condition: the
lowest order part of $g$ is completely reducible. This kind of singularity occurs e.g.
for curves ($n=2$) as $Z_{11}~jet_5(f)=z_2(z_1^3+z_2^4)$,
for surfaces ($n=3$) as $T_{455}~jet_4(f)=z_1(z_2z_3+z_1^3)$.

The lifted variety is:
\beq
\tSi(x,(l_i)_{i=1}^{q-1},\Omega^{*})=\Big\{(x,(l_i)_{i=1}^{q-1},\Omega^{*},f)\Big|\ber
f^{(p)}\sim\rm{SYM}(\Omega^{(p_1)},\dots,\Omega^{(p_k)},\Omega^{(q)})\\
\Omega^{(q)}(x)\sim\rm{SYM}(l_1,\dots,l_{q-1}),~l_i(x)=0=\Omega^{(p_i)}(x)\eer\Big\}
\eeq
\li{\bf The normals to some of the hypersurfaces are in special position.} Some simple cases:
\bei
\item{several coinciding normals: $f=\Big(\prod_{i=1}^k f_i\Big)z^r_1\Omega^{(q)},~~\Omega^{(q)}(x)\sim x^{q-1}_1$}
For example this occurs for curves ($n=2$) as $A_3,D_6,E_7,E_8,W_{13},W_{1,0},W_{17},W_{18}$,
for surfaces ($n=3$) as $X_9\dots$.

The lifted variety:
\beq
\tSi(x,l,\Omega^{*})=\Big\{(x,l,\Omega^{*},f)\Big|\ber
f^{(p)}\sim\rm{SYM}(\Omega^{(p_1)},\dots,\Omega^{(p_k)},l^r,\Omega^{(q)})\\~~
\Omega^{(q)}(x)\sim\rm{SYM}(l,\dots,l)~l(x)=0=\Omega^{(p_i)}(x)\eer\Big\}
\eeq
\li{several normals in one plane: $f=l_1\dots l_k\Omega^{(p)}$} For $n=2$ it is e.g.
an ordinary multiple point. For $n=3$ $U_{12}:~z^3_1+z^3_2+z^4_3$
\beq
\tSi(x,(l_i)_{i=1}^k,\Omega^{(p)})=\Bigg\{(x,(l_i)_{i=1}^k,f)\Big|
\ber f^{(p+k)}\sim\rm{SYM}(l_1,\dots,l_k,\Omega^{(p)})\\
(l_i(x)=0)_{i=1}^3~~~\Omega^{(p)}(x)=0~~\rm{rk}\Bigg(\begin{smallmatrix}l_1\\\dots\\l_k\end{smallmatrix}\Bigg)<3\eer\Bigg\}
\eeq
\eei
\li{\bf The intersection of two hypersurfaces lies in the singular locus of one of them}

For example: $f=z_1(z^2_2\Omega^{(p)}+z_1\Omega^{(p+1)})$.
The simplest case is $n=3:~S_{11}~jet_3(f)=z_3(z^2_2+z_1z_3)$. We omit the calculations.
\section{Some results from intersection theory}

\subsubsection{The multiplicity of intersection}\label{SecMultiplicityOfIntersection}
At each step of the degenerating process we intersect the lifted stratum $\tSi$ with a hypersurface. As
the intersection is in general non-transversal  the resulting variety will be typically reducible: in addition
to the needed (degenerated) stratum it will contain a residual variety over a cycle of jump. This residual
variety will (in general) enter with a non-trivial multiplicity. The multiplicity is calculated in the
classical way. Suppose the hypersurface is defined by the equation $\{g=0\}$. Restrict the function $g$
to the stratum $\tSi$ and find the vanishing order along the residual variety.

We illustrate this procedure in a typical:
\bex
Degeneration of the ordinary multiple point. Start from the lifted stratum:
\beq
\tSi=\{(x,f)|~~f^{(m)}|_x=0\}\subset\mP^n_x\times\mPN
\eeq
Suppose we want to degenerate by intersection with the hypersurface:
\beq
S=\{f|_x^{(p)}(\underbrace{y,\dots,y}_k)_{i_1,\dots,i_{p-k}}=0\}~~~~~m+1\le p\le l+k
\eeq
This case occurs for example in enumeration of multiple point of co-rank $r$ (in particular $A_2$ point).

The intersection $S\cap \tSi$ is non-transversal over the diagonal $x=y$. To calculate the multiplicity,
i.e. to find the
vanishing order we expand $y=x+\Delta y$. Then restricting to $\tSi$ we have (neglecting
the numerical coefficients since we are interested in the vanishing order only):
\beq\ber
\tinyM\tinyT f|_x^{(p)}(\underbrace{x+\Delta y,\dots,x+\Delta y}_k)_{i_1,\dots,i_{p-k}}
\sim
\sum_{i=0}^k f|_x^{(p)}(\underbrace{x \dots x}_{k-i}\underbrace{\Delta y\dots\Delta y}_i)_{i_1,\dots,i_{p-k}}
\sim
\tinyM\tinyT f|_x^{(p)}(\underbrace{x,\dots,x}_{p-m-1},\underbrace{\Delta y,\dots,\Delta y}_{k+1-p+m})_{i_1,\dots,i_{p-k}}+
\bet higher\\order\\terms\eet
\eer\eeq
So, the function  $f|_x^{(p)}(\underbrace{y,\dots,y}_k)_{i_1,\dots,i_{p-k}}$ has over the diagonal $x=y$ zero of the
generic order $(k+m+1-p)$. Therefore:
\beq
\tSi\cap X=\tSi_{degenerated}\cup(k+m+1-p)\tSi|_{x=y}
\eeq
\eex
\subsubsection{On cohomology classes of cycles of jump}\label{SecCohomologyClassesOfCyclesofJump}
The possible cycles of jump are described in $\S$ \ref{SecPossibleCyclesOfJump}. Here we present
their cohomology classes in the cohomology ring of their ambient space (which is the auxiliary space).
The corresponding varieties (degeneracy loci) are known classically, in particular the cohomology classes are given in
 \cite[section 14.5]{Ful}.

We define the incidence correspondence:
\beq
\Sigma=\{(y_1,\dots,y_k)|y_k\in\rm{Span}(y_1,\dots,y_{k-1})\}
\subset\mP^n_{y_1}\times\dots\times\mP^n_{y_k}
\eeq
Note, that for $k>2$ this variety is not closed. Its closure is:
\beq
\overline\Sigma=\{\mbox{The points }(y_1,\dots,y_k)\mbox{ lie in a }(k-2)\mbox{-plane}\}
\subset\mP^n_{y_1}\times\dots\times\mP^n_{y_k}
\eeq
The cohomology class of  such variety is:
\beq
[\overline\Sigma]=\sum_{i_1+\dots+i_k=n+1-k}Y_1^{i_1}Y_2^{i_2}\dots Y_k^{i_k}
\eeq
The points of $\overline\Sigma\setminus\Sigma$ correspond to configuration
$\Big\{y_1,\dots,y_{k-1}$ are linearly dependent$\Big\}$.
This subvariety (of codimension 1) will be also important in the calculations.

In case the cycle is defined by projection: $\pi_I$ (i.e. $(\{\pi_I(y^{(j)}\}_{j\in J})$ are linearly dependent),
one continues similarly (thinking of $\pi_I(y^{(j)}$ as being a point in $\pi_I(\mP^n)$).

We will often face another condition of a special type: proportionality of symmetric tensors.
Let $a,b$ be two symmetric tensors of rank $p$. By writing their independent components in a row
we can think of each of them as being a point (in homogeneous coordinates) of some big projective space $\mP^N$. Then
the proportionality of the tensors means that $a$ and $b$ coincide as the points in $\mP^N$. This condition
was considered above. Its cohomology class is:
\beq\label{EqProportionalityTwoTensors}
[a\sim b]=\sum_{i=0}^{{p+n\choose{n}}-1}A^iB^{{p+n\choose{n}}-1-i}
\eeq
Here on the right hand side $A,B$ are the cohomology classes of the elements of $a,b$ i.e. the classes of the
corresponding hypersurfaces. Equivalently, they are the first Chern classes of the corresponding line bundles.
\subsubsection{The cohomology class of a restriction of fibration}\label{SecCohomologyClassOfRestrictionFibration}
As was explained in $\S$ \ref{SecResidualVarieties} the fibration $\tSi\rightarrow Aux$ is generically
locally trivial, it is not locally trivial over the cycles of jump ($C_i\subset Aux$). The key to degenerating
procedure is the calculation of the cohomology class of the restriction $\tSi|_{C_i}\subset Aux\times\mPN$.

The
first naive idea is to represent it as a product: $[C]R_{codim\tSi-codimC}$, here $R$ is a polynomial representing a
class in $H^{2(codim\tSi-codimC)}(Aux\times\mPN)$. This happens only in special cases.
\bel
\li If $C$ is defined by a set of monomial equations in $Aux$ (i.e. by a set of the form
$\vec{x}^{\vec{m_1}}=\dots=\vec{x}^{\vec{m_k}}=0$, for $\vec{m_i}$ multi-degrees), then
$[\tSi|_C]=[C]R_{codim\tSi-codimC}$.
\li If $C$ is a "diagonal" in $Aux$ (i.e. $C_k=\{(y_1\dots y_k)|\mbox{linearly dependent}\}$) then, for a
"flag of sub-diagonals" ($C_i=\{(y_1\dots y_i)|\mbox{linearly dependent}\}$, $2\le i<k$), we have:
$[\tSi|_C]=[C_k]R_{codim\tSi-codimC}+[C_{k-1}]R_{codim\tSi-1-codimC}+\dots+
[C_2]R_{codim\tSi-n}$.
\eel
\bpr
The first claim is immediate since it follows that (up to the multiplicity) the restriction  $\tSi|_C$ lies in a linear subspace of
$Aux\times\mPN$.
To prove the second, note that $\tSi$ is defined by a collection of conditions:
$\{f^{(p)}(y_{i_1}\dots y_{i_r})_{***}\}=0$ (section \ref{SecDefinCondLinearQuasiHomog}). Thus over the open subset of $C$:
 $y_r\in span(y_1\dots y_{r-1})$ the variable $y_r$ can be eliminated from the conditions. So, the two sets of conditions
 (conditions of $C\subset Aux$ and conditions of $\tSi$ over the open subset of $C$) are explicitly transversal.
 The non-transversality can happen only over the 'infinity': $(y_1\dots y_{r-1})|\mbox{linearly dependent}$. So, we can write:
 $$[\tSi|_C]=[C_k]R_{codim\tSi-codimC}+[\mbox{a piece over }C_{k-1}].$$ Then by recursion we get the
 statement of the lemma.
\epr

An important case of the above lemma is the simplest case $C=\{y_1=y_2\}\subset Aux=\mP^n_{y_1}\times\mP^n_{y_2}$.
In this case the residual piece can be written explicitly:
\bcor\label{CorCohomologyClassResidualPieceOverDiagonal}
Suppose the projection $\tSi(y_1,y_2)|_{y_1=y_2}\mapsto\tSi(y_1)$ has the
generic fiber $\mP^r,~0\le r<n$. Then
$[\tSi(y_1,y_2)|_C]=\frac{[y_1=y_2]}{(n-r)!}\frac{\di^{n-r}[\tSi(y_1,y_2)]}{\di Y^{n-r}_2}$.
\ecor
\bpr
Over the diagonal $(y_1=y_2)$ the variable $y_2$  can be completely eliminated from the defining conditions
of $\tSi$. This corresponds to projection: $Aux=\mP^n_{y_1}\times\mP^n_{y_2}\rightarrow \mP^n_{y_1}$.
Then the class of the image is obtained by the Gysin homomorphism (section \ref{SecMinimalLifting}) from the initial class.
\epr

In general, the calculation of the cohomology class of the restriction is done as follows.
As will be
explained later we can assume $\tSi|_{C_i}$ to be irreducible (reduced). The calculation is in fact
a typical procedure from intersection theory and does not use any property of $\tSi$ related to the
singularity theory.

So, let $\tSi\subset Aux\times\mPN$ be an irreducible (reduced) projective variety. In general
$\tSi$ is not a globally complete intersection, however we assume that we can calculate the
cohomology class $[\tSi]\in H^*(Aux\times\mPN)$ by the classical intersection theory (i.e.
by intersecting various hypersurfaces and subtracting the contributions of the residual pieces). Let $\{C_i\}_i$
be all the cycles of jump.
\bel
The classes $[\tSi|_{C_i}]\in H^*(Aux\times\mPN)$ can all be calculated recursively using the following data:
\bei
\item The cohomology class of $\tSi$, obtained by the classical intersection theory
\item The cohomology classes of $C_i\subset Aux$, obtained by the classical intersection theory
\eei
\eel
\bpr
The proof goes by the induction on the grading of the cycles of jump and by the recursion on the dimensionality
of the auxiliary space $Aux$.

We first calculate $[\tSi|_{C}]$ for $C$ the cycle of jump of grading 1
(see the definition \ref{DefCycleOfJumpGrading}). As follows from proposition \ref{ClaimCycleOfJump&Hypersurf} for such
a cycle there exists a hypersurface $\{g=0\}$ containing $C$ and do not containing any other cycle. So, consider
the intersection $\tSi\cap\{g=0\}$. There can be two possibilities for the jump of dimension of fibre
 and the codimension of the cycle in $Aux=Aux_0$ (by corollary \ref{ClaimJumpOfDimensionVsCodimension}):
\bei
\item $\Delta dim_{C_i}<codim_{Aux}(C)-1.$ In this case the intersection gives just the restriction of the fibration
$\tSi\cap(g=0)=\tSi|_{g=0}$, without any residual terms.
Then decomposing the polynomial into irreducible factors $g=\prod_i g_i^{n_i}$ we get the union of restrictions:
$\tSi|_{g=0}=\bigcup_i n_i\tSi|_{g_i}$, each restriction again being irreducible.

 Consider
now the cycle $C$ as a subvariety of the new (smaller) auxiliary space $Aux_1=\{g=0\}\subset Aux_0$ intersect it with the
next hypersurface and so on.
\item $\Delta dim_{C}=codim_{Aux}(C)-1$ In this case the intersection brings residual piece (of the same dimension):
\beq\label{ClassOfRestriction}
\tSi\cap(g=0)=\tSi|_{g=0}\cup\alpha\tSi|_C
\eeq
Here $\alpha$ is the multiplicity with which
the residual piece enters. Note that the residual piece consists of the restriction to the cycle $C$ only, because
we have chosen $C$ to be of grade 1 (and then all other restrictions are excluded by codimension). In this case
we actually obtain the needed restriction as a residual piece, so the problem is reduced to the calculation of
$[\tSi|_{g=0}]$.
\eei
Thus in both cases the calculation is reduced to enumeration in a smaller auxiliary space: $Aux_1=\{g=0\}\subset Aux_0$

By the assumption of the lemma, the class $[\tSi]$ is obtained by the classical intersection procedure
(consisting of intersections and removal of residual varieties contribution). It follows that the class
$[\tSi|_{g=0}]$ can be obtained by the same procedure. In the course of calculation there will appear new
cycles of jump, however the dimension of the auxiliary space has been reduced by 1. In this way by recursion we
obtain the class of restriction: $[\tSi|_C]$ for any cycle $C$ of grade 1.

The case of general grade is treated by induction. Suppose we have calculated the classes of restriction
$[\tSi|_{C_i}]$ for all the cycles of grades up to $k$. When doing the procedure for a cycle $C$ of
grade $(k+1)$ the only difference will be that the equation (\ref{ClassOfRestriction}) is replaced by a more general:
\beq
\tSi\cap(g=0)=\tSi|_{g=0}\cup\alpha\tSi|_{C}\bigcup \alpha_i\tSi|_{C_i}
\eeq
that is on the right hand side there appear restrictions to other cycles. However (as was noted above), by lemma
\ref{ClaimCycleOfJump&Hypersurf} other cycles will be of grade at most $k$, the case already solved. So, from the
above equation we get the class of the needed restriction.
\epr
\subsubsection{On use of consistency conditions}\label{SecUseOfConsistenCondition}
In the preceding sections we have described how to calculate the cohomology classes of residual varieties.
The method is recursive and often is quite cumbersome (though it is always possible to perform the
calculations using computer). It happens, that one often can avoid lengthy calculations by using
(heavily) the consistency conditions.

The consistency conditions were stated in $\S$ \ref{SecLiftings}
(lemma \ref{ClaimFirstConsistencyCondition} and \ref{ClaimSeconConsistencCondition}). An "experimental" observation
is that they are very restrictive and in fact often {\it fix the cohomology classes of residual varieties}.
(This happens for {\bf all} the examples considered in the paper).

In general, the consistency conditions fix the cohomology class in the following equation:
\beq
[\tSi][degeneration]=[\tSi_{degenerated}]+[(y_1\dots y_k)|\mbox{linearly dependent}][R_k]+\dots+
[(y_1,y_2)|\mbox{linearly dependent}][R_2]
\eeq
Here the classes of the initial stratum ($\tSi$) and degenerating divisor/cycle are known, while the class
of degenerated stratum ($\tSi_{degenerated}$) satisfies some consistency conditions (symmetric in some variables,
with no terms of $Y^i,~~i>n-k$). The "experimental fact" ({\it which happens in all the examples considered in the paper})
is:
{\bf the above equation, together with consistency conditions has unique solution}.

The general formal way of calculations was already explained in details (sections \ref{SecPossibleCyclesOfJump},
\ref{SecIdeologyOfDegenerations}, \ref{SecCohomologyClassOfRestrictionFibration}). Thus we do not consider the
above equation
in details and do not prove any general statement of uniqueness of solution. We emphasize, however, that all
the results of this paper can be (and in fact were) obtained using the consistency conditions only.
\section{Some explicit formulae for cohomology classes of singular strata}\label{SecApendCohomologClasses}
\subsubsection{Computer calculations}\label{SecComputerCalcul}
As was already mentioned, except for the simplest cases  (ordinary multiple points, $A_{k\leq3}$) the computation
should be done on computer (systems as Singular or Mathematica can be used). The specific programs can be
obtained from the author.

Here we meet the following difficulty of purely software nature. The calculation consists of addition/subtraction
and multiplication of polynomials of indefinite degree and indefinite number of variables (i.e. both the degree
and the number of variables are parameters).  For example, when enumerating the double point of corank $r$
(section \ref{SecFormOfCorank}) the multidegree of $\tSi^n_r$ is a polynomial in $(f,x,y_1,\dots,y_r)$
of degree $n+1+r\frac{2n-1+r}{2}$. Here both $r$ and $n$ are {\it parameters}.
Another task is elimination and solution of systems of big
linear equations.
To the best of our knowledge, neither Singular nor Mathematica can in general process such expressions (i.e. open
the brackets, simplify, extract the coefficient of, say, $y^{n-r}_1\dots y^{n-r}_k$).

However the programs solve perfectly the problem for any fixed values of $n,r$. Thus to obtain the final answer
(which is a polynomial in $n,r$) one should calculate separately for a sufficient number of pairs $(n,r)$ and
then interpolate.
For the interpolation to be rigorous, one must know the degree of the polynomial we want to recreate.
This degree is known by universality \cite{Kaz2}.

We emphasize that this problem is due to the current state of software only, and not of any mathematical origin.
\subsubsection{On the possible checks of numerical results}
As in every problem, whose answer is explicit numerical formula, it is important to have some ways to
check the numerical results. Our results "successfully pass" the following checks:
\\
$\bullet$Comparison to the already known degrees. The most extensive "database" in this case is Kazarian's
tables of Thom classes for singularities of codimension up to 7. Our results are obtained
by specializations from the general case to the case of complete linear system hypersurfaces of degree $d$ in $\mP^n$.

Very few of Kazarian's results were known before.
The degree for ordinary multiple point is a classical result (known probably from the 19'th century). Another
check is for the cusp ($A_2$), enumerated by P.Aluffi.\\
$\bullet$Comparison to the known results for curves ($n=2$). By universality, the substitution $n=2$ to the formulae
must give the degrees of the strata for curves. This enables to check,
for example: $A_3,A_4,D_4,D_5,$$E_6,X_9,Z_{11},W_{12}$\\
$\bullet$Comparison to the case when the jet is reducible. For example, for singularity with degenerate quadratic form
of co-rank $k$, the substitution $k=n-2$ or $k=n-1$ gives singularities with reducible two-jet. And in these cases
the enumeration is immediate.
\subsubsection{Cohomology classes for hypersurfaces}
We present here the cohomology classes of the (minimally) lifted strata:
\beq
\tSi(x)=
\overline{\{(x,f)|~\mbox{the hypersurface defined by}~f~\mbox{has singularity of the given type at the point}~x\}}\subset
\mP^n_x\times\mPN
\eeq
The classes $[\tSi(x)]$ are expressed in terms of the generators of the cohomology ring of the ambient space
($H^*(\mP^n_x)=\mZ[X]/X^{n+1},~~H^*(\mPN)=\mZ[F]/F^{D+1}$). So, the polynomials represent the multi-degrees of the
lifted strata. The degree of the stratum itself $[\Sigma]$ is the coefficient of $X^n$.

All the notations are from \cite{AVGL}. Here (as anywhere in the paper): $d$ is the degree of singular hypersurfaces, $n$
is the dimensionality of the ambient space (thus hypersurfaces are of dimension $n-1$).
\bprop
The cohomology classes of the lifted strata and the degrees of the strata in several simplest cases are:
\eprop
$\bullet${\bf Ordinary multiple point:} $f=\sum_i z^{p+1}_i$.  Includes, for $n=2$: $A_1,D_4,X_9,\dots$; for $n=3$:
$A_1,P_8,\dots$.
$$[\tSi(x)]=Q^{n+p\choose{p}},~~Q=(d-p)X+F,~~~~~[\Sigma]={{n+p\choose{p}}\choose{n}}(d-p)^n.$$
$\bullet${\bf Degenerate multiple point (with reducible defining form):} (Defined in Appendix A)\\
$jet_p(f)=\prod_{i=1}^k\Big(\Omega^{(p_i)}_i\Big)^{r_i}$   $\sum_{i=1}^k r_ip_i=p$.
Includes:
\\$\star$ {curves $n=2$}: $A_1(k=2,r_i=1,~p_i=1),$ $A_2(k=1,r_1=2,~p_1=1),$ $D_4(k=3,r_i=1,~p_i=1),$
$E_6(k=1,r_1=3,~p_1=1),$ $X_9(k=4,r_i=1,~p_i=1),$ $Z_{11}(k=2,r_1=1,r_2=3,~p_i=1),$
$W_{12}(k=1,r_1=4,~p_1=1),\dots$
\\$\star$ {surfaces $n=3$}: $A_2(k=2,r_i=1,~p_i=1),$ $D_4(k=1,r_1=2,~p_1=1),$ $T_{3,4,4}(k=2,r_i=1,~p_1=1,p_2=2),$
$T_{4,4,4}(k=3,r_i=1,~p_i=1),$ $V_{1,0}(k=2,r_1=1,r_2=2,~p_i=1),$ $V'_1(k=1,r_1=3,~p_1=1)\dots$.

Here the forms $\Omega^{(p_i)}_i$ are mutually generic (i.e. the corresponding
hypersurfaces intersect generically near the singular point).
The cohomology (multi-)class
was given in equation (\ref{EqClassesForReducibleForms}). To obtain the answer we should extract from the equation
the coefficient of maximal (non-zero) powers of $\Omega^{(p_1)}_1\dots\Omega^{(p_k)}_k$.
\beq\hspace{-1cm}
[\tSi(x)]=\frac{1}{\rm{|Aut|}}\sum_{i=0}^{{p+n\choose{n}}-1}Q^{{p+n\choose{n}}-1-i}
X^{k+i-\sum_{j=1}^k{p_j-1+n\choose{p_j}}}
\hspace{-2.5cm}\sum_{\ber ~~~~~~~~~~~~i_1+\dots+i_k=i\\{p_j-1+n\choose{p_j-1}}-1\le i_j\le{p_j+n\choose{n}}-1\eer}
\hspace{-2.5cm}{i\choose{i_1\dots i_k}}\prod_{j=1}^k r_j^{i_j}{{p_j-1+n\choose{n}}\choose{{p_j+n\choose{n}}-1-i_j}}
\eeq
Here $Aut$ is the group of automorphisms of the branches, $i\choose{i_1~\dots~i_k}$ is the multinomial coefficient from
expansion of $(\dots)^i$ and $Q=(d-p)X+F$.
\\
$\bullet${\bf Singularity with degenerate quadratic form $\Sigma^n_k:$} $f=\sum^k_{i=1}z^3_i+\sum^n_{i=k+1}z^2_i$.
($k=1:A_2,~~k=2:D_4,~~k=3:P_8\dots$).
$$[\tSi^n_k(x)]=C_{n,k}(Q+X)^{n+1}Q^{k\choose{2}}
\sum^k_{i=0}\frac{{n-i\choose{k-i}}{k\choose{i}} }{{2k\choose{k+i}}}Q^{k-i}(-x)^i~~~~~~Q=F+(d-2)X$$\\
$$C_{n,1}=2,~~C_{n,2}=2{n+1\choose{1}},~~C_{n,3}=2{n+2\choose{3}},~~C_{n,4}=2{n+3\choose{5}}\frac{n+1}{3}$$
\\
$$C_{n,n}=2{2n\choose{n}},~~C_{n,n-1}=\frac{2^k{2k\choose{k}}}{n},~~C_{n,n-2}=\frac{{2(k+1)\choose{k+1}}{2k\choose{k}}}{{k+2\choose{2}}}$$
In particular:
\\$\star$ $[\Sigma_{A_2}]=3{n+2\choose{3}}(d-1)^{n-1}(d-2)$
\\$\star$  $[\Sigma_{D_4}]=\frac{(n+1)}{8}{n+1\choose{3}}(d-1)^{n-3}(d-2)^2\Big((d-2)(n^3+n^2+10n+8)+4(n^2+6)\Big)$
\\$\star$  $[\Sigma_{P_8}]={n+2\choose{3}}{n+4\choose{7}}(d-1)^{n-6}(d-2)^3\Big(\ber {n+7\choose{3}}\frac{(d-2)^3}{10}+{n+6\choose{2}}(d-2)^2+(n+5)\frac{9(d-2)}{2}+12\eer\Big)$
\li{\bf $A_3:$} $f=z^4_1+\sum_{i=2}^nz^2_i,~~~~Q=F+(d-2)X,~~~~
[\tSi_{A_3}(x)]=(Q+X)^{n+1}(nQ-2X)(Q\frac{3n-1}{2}-4X)$.
$$[\Sigma_{A_3}]={n+2\choose{3}}(d-1)^{n-2}\Bigg(\frac{(3n-1)(n+3)}{2}(d-2)^2+2(n-1)(d-2)-4\Bigg)$$
$\bullet${\bf $A_4:$} $f=z^5_1+\sum_{i=2}^nz^2_i$.
$$[\tSi_{A_4}(x)]=(Q+X)^{n+1}\Big(\frac{n(5n^2-5n+2)}{2}Q^3-4(5n^2-3n+1)Q^2x+4(13n-5)Qx^2-48x^3\Big)~~~~Q=F+(d-2)X$$
$$[\Sigma_{A_4}]=\frac{1}{8}{n+2\choose{3}}(d-1)^{n-3}
\Bigg(\ber d^3(24-46n+27n^2+30n^3+5n^4)-d^2(+96-184n+138n^2+160n^3+30n^4)+\\+d(144-316n+192n^2+280n^3+60n^4)
-136+224n-48n^2-160n^3-40n^4\eer\Bigg)$$
$\bullet${$D_5$:} $f=z^4_1+z_1z^2_2+\sum_{i=3}^nz^2_i,~~~~Q=(d-2)X+F$.
$$[\tSi_{D_5}(x)]=(Q+X)^{n+1}\frac{Q(n+1)}{6}\Bigg((3n-2)Q-10X\Bigg)\Bigg((n^2-n)Q^2-6(n-1)QX+12X^2\Bigg)$$
$\bullet${$D_6$:} $f=z^5_1+z_1z^2_2+\sum_{i=3}^nz^2_i,~~~~Q=(d-2)X+F$
$$[\tSi_{D_6}(x)]=(Q+X)^{n+1}(1+n)Q
\Bigg(\ber\frac{4(n-1)n(3n^2-5n+3)}{15}Q^4-\frac{2(n-1)(16n^2-19n+9)}{3}Q^3x+\\
\frac{2(83n^2-140n+69)}{3}Q^2x^2-136(n-1)Qx^3+136x^4\eer\Bigg)$$
$\bullet${\bf $E_6$:} $f=z^4_1+z^3_2+\sum_{i=3}^nz^2_i,~~~~Q=(d-2)X+F$.
$$[\tSi_{E_6}(x)]=(Q+X)^{n+1}Q(n+1)\Bigg(\ber\frac{n(n-1)(12n^2-15n+2)}{40}Q^4-\frac{(n-1)n(15n-14)}{4}Q^3X+\\+
\frac{37n^2-52n+12}{2}Q^2X^2-6(7n-5)QX^3+36X^4\eer\Bigg)$$
$\bullet${\bf $X_9$:} $f=z^4_1+z^4_2+\sum_{i=3}^nz^2_i,~~~~Q=(d-2)X+F$.
$$[\tSi_{X_9}(x)]=(Q+X)^{n+1}Q(n+1)\Bigg(\ber\frac{n(n-1)(3n+1)(33n^3-102n^2+102n-20)}{1680}Q^6-
\frac{(n-1)n(138n^3-309n^2+153n+86)}{120}Q^5X\\+
\frac{(379n^4-1004n^3+719n^2+166n-160)}{40}Q^4X^2-\frac{126n^3-247n^2+58n+72}{2}Q^3X^3+\\+
\frac{625n^2-670n-216}{6}Q^2x^4-16(8n-1)QX^5+48X^6\eer\!\!\!\!\Bigg)\\
$$
$\bullet${\bf $Q_{10}$:} $z_1^4+z_2^3+z_1z_3^2+\sum_{i=4}^n z_i^2,~~~~Q=(d-2)X+F$.
$$[\tSi_{Q_{10}}(x)]=(Q+X)^{n+1}{n+2\choose{3}}Q^3\Big((n-1)Q-4X\Big)^2\Bigg(\frac{3}{5}{n\choose{3}}Q^3-
\frac{12}{5}{n-1\choose{2}}Q^2X+6(n-2)QX^2-12X^3\Bigg)$$
$\bullet${\bf $S_{11}$:} $z_1^4+z_2^2z_3+z_1z_3^2+\sum_{i=4}^n z_i^2,~~~~Q=(d-2)X+F$.
$$[\tSi_{S_{11}}(x)]=(Q+X)^{n+1}6{n+2\choose{3}}Q^3\Big((n-1)Q-4x\Big)
\Bigg(\ber\frac{(n-2)(n-1)n(51n^2-98n+31)}{3360}Q^5+\frac{(56n-79)}{2}Qx^4
\\-\frac{(n-2)(n-1)(253n^2-392n+75)}{840}Q^4x-27x^5
\\+\frac{(n-2)(103n^2-190n+67)}{40}Q^3x^2-\frac{(47n^2-130n+75)}{4}Q^2x^3\eer\Bigg)$$
$\bullet${\bf $U_{12}$:} $z_1^3+z_2^3+z_3^4+\sum_{i=4}^n z_i^2,~~~~Q=(d-2)X+F$.
$$[\tSi_{U_{12}}(x)]=(Q+X)^{n+1}{n+2\choose{3}}Q^3
\Bigg(\ber
{n\choose{3}}\frac{(n-1)(117n^3-328n^2+207n+12)}{1120}Q^7-
{n-1\choose{2}}\frac{(323n^4-1057n^3+1021n^2-231n-72)}{336}Q^6x\\
+{n-1\choose{2}}\frac{(991n^3-2545n^2+1525n-27)}{84}Q^5x^2+
\frac{(n-1)(2161n^2-5606n+2553)}{12}Q^3x^4
\\
-\frac{(3491n^4-16290n^3+25396n^2-14838n+2577)}{84}Q^4x^3+28(25n-31)Qx^6-\\
(475n^2-1150n+631)Q^2x^5-448x^7
\eer
\Bigg)$$

{\it Address}: School of Mathematical Sciences, \mbox{Tel Aviv}
University, \mbox{Ramat Aviv}, 69978 \mbox{Tel Aviv}, Israel.

{\it E-mail}: kernerdm@post.tau.ac.il

\begin{thebibliography}{99}\small\addtolength{\parskip}{-8pt}
\bibitem[AGLV]{AVGL} V.I.Arnol'd, V.A.Vasil'ev, V.V.Goryunov, O.V.Lyashko,  {\it Singularities. I. Local
and global theory.} Current problems in mathematics. Fundamental directions, Vol. 6, 5--257, Itogi Nauki i Tekhniki,
Akad. Nauk SSSR, Vsesoyuz. Inst. Nauchn. i Tekhn. Inform., Moscow, 1988.
\bibitem[Aluf98]{Alufi}P.Aluffi {\it Characteristic classes of discriminants and enumerative geometry}, Comm. Algebra
26 (1998), no. 10, 3165--3193.
\bibitem[Dimca]{Dim} A.Dimca {\it Singularities and topology of hypersurfaces}, Universitext. Springer-Verlag, New York, 1992.
\bibitem[EyrGas05]{EG} C.Eyral, E.Gasparim {\it Multiplicity of complex hypersurface singularities, Rouche' satellites
and Zariski's problem}, preprint math.AG/0509409
\bibitem[Fult]{Ful} W.Fulton, {\it Intersection theory}, Second edition.
 A Series of Modern Surveys in Mathematics, Springer-Verlag, Berlin, 1998.
\bibitem[Greue86]{Greuel} G.-M.Greuel, {\it Constant Milnor number implies constant multiplicity for quasihomogeneous
singularities}. Manuscripta Math. {\bf 56 (1986)}, no. 2, 159--166.
\bibitem[GreLoShu06]{Shustin} G.-M.Greuel, C.Lossen, E.Shustin {\it Introduction to Singularities and Deformations}.
Series: Springer Monographs in Mathematics 2006. ISBN: 3-540-28380-3
\bibitem[GreLoShu01]{GLS} G.-M.Greuel, C.Lossen, E.Shustin {\it The variety of plane curves with ordinary singularities
is not irreducible} Internat. Math. Res. Notices {\bf 11} (2001), 543--550.
\bibitem[GrePfi96]{GreuelPfister1}G.-M.Greuel, G.Pfister {\it Advances and improvements in the theory of
standard bases and syzygies} Arch. Math.{\bf 66 (1996)}, 163-176.
\bibitem[HerV\'az01]{HV}R.Hern\'andez M.J.V\'azquez-Gallo {\it Degree of strata of singular cubic surfaces}
Trans. Amer. Math. Soc. 353 (2001), no. 1, 95--115.
\bibitem[Kaz00]{Kaz1} M.Kazarian {\it Classifying spaces of singularities and Thom polynomials} in "New developments in
Singularity Theory (Cambridge 2000)", NATO Sci. Ser. II Math. Phys.
Chem, 21, Kluwer Acad. Publ.,  Dordrecht, 2001, pp. 117--134.
\bibitem[Kaz03-1]{Kaz2} M.Kazarian {\it Thom polynomials for Lagrange, Legendre, and critical point function
singularities}, Proc. Lond. Math. Soc. (3) {\bf 86} (2003), 707--734.
\bibitem[Kaz03-2]{Kaz3} M.Kazarian {\it Multisingularities, cobordisms, and enumerative geometry},
 Russ. Math. Surveys {\bf 58(4)} (2003), 665--724.
\bibitem[Kaz]{Kaz4} M.Kazarian {\it Characteristic Classes in Singularity theory},
Doctoral Dissertation (habilitation thesis), {\it Steklov Math. Inst.},2003.
\bibitem[Ker06]{Ker}D.Kerner {\it Enumeration of singular algebraic curves}, Israel
Journal of Math. 155 (2006), pp1-56,  arXiv:math.AG/0407358
\bibitem[Klei76]{Klei1} S.L.Kleiman, {\it The enumerative theory of singularities}. Real and complex singularities
(Proc. Ninth Nordic Summer School/NAVF Sympos. Math., Oslo, 1976), pp. 297--396. Sijthoff and Noordhoff,
Alphen aan den Rijn, 1977.
\bibitem[Klei87]{Klei2} S.L.Kleiman {\it Intersection theory and enumerative geometry: a decade in review}.
With the collaboration of Anders Thorup on §3. Proc. Sympos. Pure Math., 46, Part 2, Algebraic geometry, Bowdoin, 1985
(Brunswick, Maine, 1985), 321--370, Amer. Math. Soc., Providence, RI, 1987.
\bibitem[KleiPie98]{KleiPien1} S.Kleiman, R.Piene {\it Enumerating singular curves on surfaces},
 Algebraic geometry: Hirzebruch 70 (Warsaw, 1998), 209--238, Contemp. Math., 241, Amer. Math. Soc.,Providence, RI, 1999.)
\bibitem[Lue87]{Luengo} I.Luengo {\it The $\mu$-constant stratum is not smooth.} Invent. Math. {\bf 90} (1987), no. 1, 139--152
\bibitem[Pham70]{Pham}F.Pham {\it Remarque sur l'\'{e}quisingularit\'{e} universelle}, Facult\'{e} des Sciences de Nice, 1970 (Preprint).

\bibitem[Saito71]{Saito71}K.Saito, {\it Quasihomogene isolierte Singularitäten von Hyperflächen}.
Invent. Math. 14 (1971), 123--142
\bibitem[Thom54]{Thom54}R.Thom, {\it Quelques propriétés globales des variétés différentiables.}Comment. Math. Helv. 28, (1954). 17--86

\bibitem[Vain03]{Vain2}I.Vainsencher {\it Hypersurfaces with up to six double points}.
 Special issue in honor of Steven L. Kleiman. Comm. Algebra 31(2003), no. 8, 4107--4129.
\bibitem[Var82]{Var} A.N.Varchenko {\it A lower bound for the codimension of the stratum $\mu\!=$const in terms of the mixed
Hodge structure}. Moscow University Mathematics Bulletin 37 (1982), no. 2, 28-31.
\end{thebibliography}
\end{document}